\documentclass[12pt,reqno]{amsart}

\usepackage{amssymb}
\usepackage{amsmath}
\usepackage{amsthm}

\newcommand{\nc}{\newcommand}
\nc{\nen}{\newenvironment}

\nc{\notebox}[1]{\noindent\fbox{\parbox{12.5cm}{\sf #1}}\\[8pt]}
\nc{\Thm}[1]{Theorem~\ref{#1}}
\nc{\Prop}[1]{Proposition~\ref{#1}}
\nc{\Lem}[1]{Lemma~\ref{#1}}
\nc{\Cor}[1]{Corollary~\ref{#1}}
\nc{\Conj}[1]{Conjecture~\ref{#1}}
\nc{\Claim}[1]{Claim~\ref{#1}}
\nc{\Defn}[1]{Definition~\ref{#1}}
\nc{\Exa}[1]{Example~\ref{#1}}
\nc{\Rem}[1]{Remark~\ref{#1}}
\nc{\Note}[1]{Note~\ref{#1}}
\nc{\bref}[1]{(\ref{#1})}

\nen{thm}[1]{\label{#1}{\bf Theorem.\ } \em}{}
\nen{propn}[1]{\label{#1}{\bf Proposition.\ } \em}{}
\nen{lemma}[1]{\label{#1}{\bf Lemma.\ } \em}{}
\nen{klem}[1]{\label{#1}{\bf Key Lemma.\ } \em}{}
\nen{cor}[1]{\label{#1}{\bf Corollary.\ } \em}{}
\nen{conj}[1]{\label{#1}{\bf Conjecture.\ } \em}{}
\nen{claim}[1]{\label{#1}{\bf Claim.\ } \em}{}

\nen{defn}[1]{\label{#1}{\bf Definition.\ } }{}
\nen{exa}[1]{\label{#1}{\bf Example.\ } }{}

\nen{remark}[1]{\label{#1}{\em Remark.\ } }{}
\nen{note}[1]{\label{#1}{\em Note.\ } }{}
\nen{exer}[1]{\label{#1}{\em Exercise.\ } }{}

\renewcommand{\frak}{\mathfrak}

\renewcommand{\(}{\left(}
\renewcommand{\)}{\right)}
\newcommand{\g}{\mathfrak{g}}
\newcommand{\hh}{\mathfrak{h}}
\newcommand{\mm}{\mathfrak{m}}
\newcommand{\lf}{\mathfrak l}

\renewcommand{\aa}{\alpha}
\newcommand{\ff}{\varphi}

\newcommand{\Uh}{{U}_h (\g)}
\newcommand{\Ug}{U(\g)}

\newcommand{\A}{\mathcal{A}}

\newcommand{\C}{\mathbb{C}}
\renewcommand{\SS}{\mathbb{S}}
\newcommand{\CC}{\mathcal{C}}

\newcommand{\EE}{\mathcal{E}}

\newcommand{\R}{\mathbb{R}}

\newcommand{\HH}{\mathcal{H}}
\renewcommand{\L}{\mathcal{L}}

\newcommand{\ZZ}{\mathbb{Z}}

\newcommand{\ot}{\otimes}

\renewcommand{\[}{[\![}
\renewcommand{\]}{]\!]}
\newcommand{\half}{\frac{1}{2}}
\newcommand{\hlf}[1]{\frac{#1}{2}}
\newcommand{\third}{\frac{1}{3}}

\newcommand{\ad}{\operatorname{ad}}
\newcommand{\Ad}{\operatorname{Ad}}

\renewcommand{\Im}{\operatorname{Im}}
\newcommand{\Ker}{\operatorname{Ker}}

\renewcommand{\ll}{\lambda}

\newcommand{\Ga}{\Gamma}

\newcommand{\Ob}{\overline\Omega}
\newcommand{\Pib}{{\overline\Pi}}

\newcommand{\alb}{{\bar\alpha}}
\newcommand{\beb}{{\bar\beta}}
\newcommand{\Vect}{\operatorname{Vect}}

\newcommand{\be}[1]{\begin{eqnarray#1}}
\newcommand{\ee}[1]{\end{eqnarray#1}}

\newcommand{\tor}[1]{\stackrel{#1}{\longrightarrow}}

\newcommand{\De}{\Delta}
\newcommand{\de}{\delta}

\renewcommand{\t}{\otimes}

\newcommand{\bbe}{{\bar\beta}}
\newcommand{\bga}{{\bar\gamma}}
\newcommand{\Ea}{E_{\aa}}
\newcommand{\Ema}{E_{-\aa}}

\newcommand{\baa}{{\bar\aa}}

\numberwithin{equation}{section}

\newcommand{\pa}{\subsection{}}
\newcommand{\pat}[1]{\subsection{\bf #1}}
\newcommand{\spa}{\subsubsection{}}
\newcommand{\spat}[1]{\subsubsection{\bf #1}}

\begin{document}

\title[$\Uh$ invariant Quantization]{$\Uh$ invariant Quantization of Coadjoint Orbits and
Vector Bundles over them}
                          
\author[Joseph Donin]{Joseph Donin\\
\bigskip
{\Tiny Dept. of Math.  Bar-Ilan University\\
Max-Planck-Institut f\"ur Mathematik}}
\address{Dept. of Math.  Bar-Ilan University, 52900 Ramat-Gan, Israel}
\curraddr{Max-Planck-Institut f\"ur Mathematik, Vivatsgasse 7, 53111 Bonn}
\email{donin@mpim-bonn.mpg.de\\
\
donin@macs.biu.ac.il}

\thanks{This research is partially supported
by Israel Academy of Sciences grant no. 8007/99-01}

\begin{abstract}
Let $M$ be a coadjoint semisimple orbit of a simple Lie group $G$.
Let $U_h(\g)$ be a quantum group corresponding to $G$.
We construct a universal family of
$U_h(\g)$ invariant quantizations of the sheaf of
functions on $M$ and describe all such quantizations.
We also describe all two parameter $U_h(\g)$ invariant quantizations
on $M$, which can be considered as $U_h(\g)$ invariant
quantizations of the Kirillov-Kostant-Souriau (KKS) Poisson bracket on $M$.
We also consider how those quantizations relate to
the natural polarizations of $M$ with respect to the KKS bracket.
Using polarizations, we quantize the sheaves of sections of 
vector bundles on $M$ as one- and two-sided $U_h(\g)$ invariant 
modules over a quantized function sheaf.
\end{abstract} 

\maketitle
\section{Introduction}

Let $G$ be a simple Lie group with Lie algebra $\g$, $M$ a semisimple coadjoint orbit of $G$,
i.e., the orbit of $G$ passing through a semisimple element in the coadjoint
representation $\g^*$.
Let $\A$ be the sheaf of functions on $M$. It may be the sheaf of smooth,
analytic, or algebraic functions. 
The universal enveloping algebra $\Ug$ acts on the sections of $\A$ and
the multiplication in $\A$ is $\Ug$ invariant.
Let $\Uh$ denote a quantum group
that is a deformation of the $\Ug$ as a bialgebra. 
In the paper \cite{DGS} we considered the following problems.

1) Does there exists a $\Uh$ invariant deformation quantization of $\A$,
i.e., a quantization, $\A_h$, having a $\Uh$ invariant multiplication?

2) Does there exists a two parameter (double) $\Uh$ invariant quantization, $\A_{t,h}$,
such that $\A_{t,0}$ is a $\Ug$ invariant quantization of $\A$ with
Poisson bracket being the Kirillov-Kostant-Souriau (KKS)
Poisson bracket on $M$? Note that $\A_{t,h}$ 
can be considered as a $U_h(\g)$ invariant
quantization of the KKS Poisson bracket on $M$.

In \cite{DGS}, we have classified the Poisson brackets admissible for
one and two parameter quantizations. 
An admissible Poisson bracket for one parameter
quantization is the same as a Poisson bracket making $M$
into a Poisson manifold with Poisson action of $G$, where
$G$ is considered to be the Poisson-Lie group with
Poisson structure defined by the $r$-matrix related to $\Uh$.
An admissible Poisson bracket for two parameter quantization,
in addition, must be compatible with the KKS bracket on $M$.

We have shown that all
semisimple orbits have admissible Poisson brackets for
one parameter quantization and almost all such brackets
can be quantized. 

In \cite{DGS}, we called an orbit having a Poisson bracket
admissible for a two parameter quantization {\em a good orbit}.
We have classified the good semisimple orbits for all simple $\g$
and shown that in case $\g\neq sl(n)$ all the good orbits can be
quantized. In the case $\g=sl(n)$ all semisimple orbits are good
but in \cite{DGS} we did not prove the existence of their
double quantization.

In this paper we give a complete description of one and two parameter
$\Uh$ invariant quantizations on $M$.
We show that for each semisimple orbit $M$ there exists a universal
family of quantizations of $\A$. This family is given by a family
of multiplications
\be{*}
m_{f,h}:\A\ot\A\to\A[[h]], \qquad f\in X,
\ee{*}
where $X$ is the manifold of all admissible Poisson brackets on $M$.
The universality means that any one parameter $\Uh$ invariant
quantization of $\A$ is given by the multiplication of the form
$m_{f(h),h}$, where $f(h)$ is a formal path in $X$, and two
different paths give nonequivalent quantizations.
As a consequence we obtain that any admissible Poisson bracket on $M$
can be quantized. 

There is the analogous description for two parameter quantizations
on any good orbit. In particular, we prove that any good
orbit (including all semisimple orbits in $sl(n)^*$) admits
a two parameter $\Uh$ invariant quantization.  

Further, we consider the natural polarizations on $M$
with respect to the KKS Poisson bracket.
We show that all $\Uh$ invariant quantizations on $M$
being restricted to functions constant along a polarization
have a standard form.
In some cases, when $M$ is a coadjoint orbit of a real Lie group $G$,
the polarizations define complex structures on $M$.
In such cases the sheaf of functions constant along polarization
specializes to the sheaf of holomorphic (or antiholomorphic)
functions on $M$. So we obtain that in the real case any quantization of
smooth functions on $M$ induces a unique quantization
of the subsheaves of holomorphic and antiholomorphic functions on $M$.

In the paper we also consider the quantization of $G$ invariant
vector bundles. We identify such a bundle with the sheaf of its
smooth sections. Let $\A_h$ be a $\Uh$ invariant quantization of
the sheaf of smooth functions on $M$. Under the quantization of 
a vector bundle $V$ on $M$ with respect to $\A_h$ we mean the sheaf
$V[[h]]$ endowed with a structure of $\Uh$ invariant left
(right, two-sided) $\A_h$ module.
Using a complex polarization, we show that for any $\Uh$ invariant quantization $\A_h$ and
any $G$ invariant vector bundle $V$ there exists a $\Uh$ invariant
quantization of $V$ as a left $\A_h$ module.
Moreover, we show that there exists a special quantization, $\A_h$,
such that any vector bundle $V$ admits a quantization as
a two-sided module with respect to $\A_h$.
Note that  the
papers \cite{Jo}, \cite{DG}, where the algebra of global holomorphic sections of
linear vector bundles on flag varieties was quantized, relate 
to the problem of quantizing vector bundles. 

The paper is organized as follows. In Section 2 we recall some
facts on quantum groups essential for our approach to $\Uh$ invariant
quantization. In particular, we define a quantum group, $\Uh$,
for any classical $r$-matrix $r$ and show that the problem of
constructing $\Uh$ invariant quantization is equivalent to
the problem of constructing $\Ug$ invariant $\Phi_h$ associative
quantization, where $\Phi_h\in\Ug^{\ot 3}[[h]]$ defines the Drinfeld
associativity constraint (see \cite{Dr2}).
Thus, a $\Ug$ invariant $\Phi_h$ associative quantization
defines a $\Uh$ invariant quantization for all quantum groups
associated with different $r$. In this section we also define 
$\ff$-brackets, which are infinitesimal parts 
of $\Ug$ invariant $\Phi_h$ associative quantizations. 
We show that a Poisson bracket admissible for $\Uh$ invariant
quantization is the difference of a $\ff$-bracket and the 
bracket induced on $M$ by the $r$-matrix associated with $\Uh$.

In Section 3 we give a classification of $\ff$- and good brackets
on semisimple orbits. In particular, we give a description of the variety of
those brackets, which is more detailed than in \cite{DGS}. We are needed
in this description in the following sections.

In Section 4 we consider Poisson cohomologies of some
parameterized complexes, which we use in Sections 5 and 6 for proving the existence
of the universal one and two parameter quantizations.

In Sections 7 we consider the natural polarizations on $M$ and
prove that any $\Ug$ invariant $\Phi_h$ associative quantization
is trivial when restricted to the sheaf of functions constant along
polarization.

Note that up to Section 7 we assume that $G$ is a complex Lie group,
$M$ is a complex subvariety in $g^*$, and $\A$ is the sheaf of complex
analytic functions on $M$.

In Section 8 we specify our results to a real Lie group $G$.
We consider complex structures corresponding to polarizations
and use them to construct the quantizations of $G$ invariant
vector bundles.

{\bf Acknowledgments.} I am very grateful to P.Bressler and D.Gurevich for
useful discussions. I thank Max-Planck-Institut f\"ur Mathematik for
hospitality and very stimulating working atmosphere.

\section{Preliminaries}

\pat{Quantum groups.}
We will consider quantum groups in sense of Drinfeld, \cite{Dr2},
as deformed universal enveloping algebras.  
If $U(\g)$ is the universal enveloping algebra of a complex Lie algebra  $\g$, then
the quantum group (or quantized universal enveloping algebra) corresponding
to $U(\g)$ is a topological Hopf algebra, $\Uh$, over $\C[[h]]$, isomorphic to $U(\g)[[h]]$ as a
topological $\C[[h]]$ module and such that $\Uh/h\Uh=U(\g)$ as a Hopf algebra
over $\C$.
In particular, the deformed comultiplication in  $\Uh$ has the form
\be{}
\label{comult}
\De_h=\De+h\De_1+o(h),
\ee{}
where $\De$ is the comultiplication in the universal enveloping algebra $U(\g)$.
One can prove, \cite{Dr2},  that the map
$\De_1:U(\g)\to U(\g)\ot U(\g)$ is such that $\De_1-\sigma\De_1=\de$
($\sigma$ is the usual permutation) being restricted to $\g$ gives
a map $\de:\g\to\wedge^2\g$ which is a 1-cocycle and defines
the structure of a Lie coalgebra on $\g$ (the structure of a Lie algebra on
the dual space $\g^*$).  
The pair $(\g,\de)$ is called
a quasiclassical limit of $\Uh$.

In general, a pair $(\g,\de)$, where
$\g$ is a Lie algebra and $\de$ is such a 1-cocycle, is called a Lie bialgebra.
It is proven, \cite{EK}, that any Lie bialgebra $(\g,\de)$ can be quantized, i.e., 
there exists a quantum group  $\Uh$ such
that  the  pair  $(\g,\de)$ is its quasiclassical limit.

A Lie bialgebra $(\g,\de)$ is said to be  a coboundary one if  there exists an element
$r\in\wedge^2\g$, called the classical $r$-matrix, 
such that $\de(x)=[r,\De(x)]$ for $x\in\g$. Since $\de$ defines a Lie coalgebra
structure, 
$r$ has to satisfy the so-called classical Yang-Baxter equation
which can be written in the form
\be{}
\label{cYB}
\[r,r\]=\ff,
\ee{}
where $\[\cdot,\cdot\]$ stands for the Schouten bracket and 
$\ff\in\wedge^3\g$ is an invariant element. We denote the coboundary
Lie bialgebra by $(\g,r)$.

In case $\g$ is a simple Lie algebra, the most known $r$-matrix is
the Sklyanin-Drinfeld one:
\be{*}
r=\sum_\alpha X_\alpha\wedge X_{-\alpha},
\ee{*} 
where the sum runs over all positive roots; the root vectors $X_\alpha$
are chosen is such a way that $(X_\alpha, X_{-\alpha})=1$
for the Killing form $(\cdot,\cdot)$.
This is the only $r$-matrix of weight zero, \cite{SS}, and its quantization is the Drinfeld-Jimbo
quantum group.
A classification of all $r$-matrices for simple Lie algebras was
given in \cite{BD}.

From results of Drinfeld and of Etingof and Kazhdan one can derive
the following
\spa\begin{propn}{propo2.1}
Let $\g$ be a semisimple Lie algebra. Then 

a) any Lie bialgebra $(\g,\de)$ is a coboundary one;

b) the quantization, $\Uh$, of any coboundary Lie bialgebra $(\g,r)$ exists and
is isomorphic to $\Ug[[h]]$ as a topological $\C[[h]]$ algebra;

c) the comultiplication in $\Uh$ has the form
\be{}
\label{comul}
\De_h(x)=F_h\De(x)F^{-1}_h, \qquad x\in\Ug,
\ee{}
 where $F_h\in \Ug^{\ot 2}[[h]]$ and can be chosen in the form
\be{}
\label{F}
F_h=1\ot 1+\hlf{h}r +o(h).
\ee{} 
\end{propn}

\begin{proof}
a) follows from the fact that $H^1(\g,\wedge^2\g)=0$.
It follows from $H^2(\g,\Ug)=0$ that
$\Ug$ does not admit any
nontrivial deformations as an algebra, (see \cite{Dr1}), which
proves b). 
From the fact that $H^1(\g,\Ug^{\ot 2})=0$ it follows that
any deformation of the algebra morphism
$\De:\Ug\to\Ug\ot\Ug$ appears as a conjugation of $\De$.
In particular, the comultiplication in $\Uh$ looks like (\ref{comul}) 
with some $F_h$ such that $F_0=1\ot 1$.
It follows from the coassociativity of $\De_h$ that $F_h$  
satisfies  the
equation
\be{}
\label{eqF}
(F_h\t 1)\cdot(\De\t id)(F_h)=(1\t F_h)\cdot(id\t\De)(F_h)\cdot\Phi_h
\ee{}
for some invariant element $\Phi_h\in\Ug^{\ot 3}[[h]]$.

The element $F_h$ satisfying (\ref{comul}) and (\ref{F}) can be obtained by a correction
of some $F_h$ only obeying (\ref{eqF}), \cite{Dr2}.
This procedure also uses a simple cohomological argument, which 
proves c).
\end{proof}

\spa It follows from (\ref{eqF}) that if $F_h$ has the form (\ref{F}), then
the coefficient by $h$ in $\Phi_h$ vanishes. Moreover, the coefficient by $h^2$
is the element $\ff$ from (\ref{cYB}), i.e.,  
\be{}\label{fPhi}
\Phi_h=1\t 1\t 1+h^2\ff+o(h^2).
\ee{}
In addition, it follows  from (\ref{eqF}) that
$\Phi_h$ satisfies the pentagon identity
\be{*}
(id^{\t 2}\t \De)(\Phi_h)\cdot(\De\t id^{\t 2})(\Phi_h)=
(1\t\Phi_h)\cdot(id\t\De\t id)(\Phi_h)\cdot(\Phi_h\t 1).
\ee{*}

\pat{Equivariant deformation quantization.}
Let $G$ be a simple connected complex Lie group whose Lie algebra is $\g$.
Let $G$ act on a manifold $M$ and $\A$ be the sheaf of 
functions on $M$. It may be the sheaf of analytic, smooth, or  
algebraic functions, dependingly of the type of $M$.
Then $\Ug$ acts on sections of $\A$, and the multiplication
in $\A$ is $\Ug$ invariant.

The deformation quantization of $\A$ is a sheaf of associative
algebras, $\A_h$, which is 
isomorphic to $\A[[h]]=\A\t\C[[h]]$ (completed tensor
product) as a $\C[[h]]$-module, with multiplication in
$\A_h$ having  the form
$m_h=\sum_{k=0}^\infty h^k m_k$, where
$m_0$ is the usual commutative multiplication in $\A$
and $m_k$, $k>0$, are bidifferential operators vanishing on constants.
The algebra $\Ug[[h]]$ is clearly acts on the $\C[[h]]$ module $\A_h$. 

We will study  quantizations of $\A$ which are invariant under
the $\Uh$ action, i.e., under the comultiplication $\De_h$. This means that
\be{*}
b m_h(x\otimes y)=m_h\De_h(b)(x\t y)\quad \mbox{for}\ b\in\Ug,\ x,y\in \A.
\ee{*}

A $\C[[h]]$ linear map $\mu_h:\A_h\ot\A_h\to\A_h$ is called
a $\Phi_h$ {\em associative} multiplication if
\be{*}
\mu_h(\Phi_1x\otimes \mu_h(\Phi_2 y\t \Phi_3 z)))=\mu_h(\mu_h(x\t y)\t z)
\quad \mbox{ for } x,y,z\in \A,  
\ee{*}
where $\Phi_h=\Phi_1\ot\Phi_2\ot\Phi_3$ (summation implicit).

We say that the $\Phi_h$ associative multiplication 
$\mu_h=\sum_{k=0}^\infty h^k \mu_k$ gives a $\Phi_h$ associative
quantization of $\A$ if $\mu_0=m_0$,  the usual multiplication in $\A$,
and $\mu_k$, $k>0$, are bidifferential operators vanishing on constants.

\spa\begin{propn}{propo2.2} There is a natural one-to-one correspondence between
$\Uh$ invariant and $\Ug$ invariant $\Phi_h$ associative quantizations
of $\A$. Namely, if $\mu_h$ is a $\Ug$ invariant $\Phi_h$
associative multiplication in $\A[[h]]$, then
\be{}\label{cormul}
m_h=\mu_hF_h^{-1}
\ee{} 
gives a $\Uh$ invariant associative
multiplication in $\A[[h]]$.    
\end{propn}
\begin{proof} This follows immediately from \bref{comul} and \bref{eqF}.
This follows also from the categorical interpretation of
$\Phi_h$ and $F_h$, \cite{Dr2}, \cite{DGS}.
\end{proof}

This proposition shows that given a $\Ug$ invariant $\Phi_h$
associative quantization of $\A$, we can get the $\Uh$ 
invariant quantization of $\A$ for any quantum
group $\Uh$ from \Prop{propo2.1} b) by applying
$F_h$ from \bref{F} to the $\Phi_h$ associative multiplication. 

\pat{Poisson brackets associated with equivariant quantizations.}
A skew-symmetric map $f:\A^{\ot 2}\to\A$ we call a
{\em bracket} if it satisfies the Leibniz rule:
$f(ab,c)=af(b,c)+f(a,c)b$ for $a,b,c\in\A$.
It is easy to see that any bracket is presented by a bivector
field on $M$. Further we will identify brackets and
bivector fields on $M$.

For an element $\psi\in\wedge^k\g$ we denote by 
$\psi_M$ the $k$-vector field on $M$ which is induced
by the action map $\g\to\Vect(M)$.

A bracket $f$ is a Poisson one if the Schouten bracket
$\[f,f\]$ is equal to zero.

\spa\begin{defn}{defphi} A $G$ invariant bracket $f$ on $M$ we call a
$\ff$-{\em bracket} if $\[f,f\]=-\ff_M$, where
$\ff\in\wedge^3\g$ is an invariant element.
\end{defn}

\spa\begin{propn}{propo2.3}
Let $\A_h$ be a $\Ug$ invariant $\Phi_h$ associative
quantization with multiplication 
$\mu_h=m_0+h\mu_1+o(h)$, where $m_0$ is
the multiplication in $\A$.
Then the map $f:\A^{\ot 2}\to\A$, $f(a,b)=\mu_1(a,b)-\mu_1(b,a)$,
is a $\ff$-bracket for $\ff$ from \bref{fPhi}.
\end{propn}

\begin{proof} A direct computation.
Another proof is found in \cite{DGS}.
\end{proof}

\spa\begin{cor}{coro2.3} Let $\A_h$ be a $\Uh$ invariant associative
quantization with multiplication $m_h=m_0+h m_1+o(h)$.
Then the corresponding Poisson bracket $p(a,b)=m_1(a,b)-m_1(b,a)$
has the form
\be{}\label{dopq}
p(a,b)=f(a,b)-r_M(a,b),
\ee{}
where $r$ is the $r$-matrix
corresponding to $\Uh$ and $f$ is a $\ff$-bracket
with $\ff=\[r,r\]$.
\end{cor}

\begin{proof} By \Prop{propo2.2} there is a $\Ug$ invariant
$\Phi_h$ associative multiplication $\mu_h$ such that
$m_h=\mu_hF_h^{-1}$ with $F_h$ as in \bref{F}.
Let $f$ be the $\ff$-bracket corresponding to $\mu_h$.
Then a direct computation shows that the Poisson bracket of
$m_h$ is as required.
\end{proof}

\spat{\em Remark.}\label{rem2.4} For an $r$-matrix $r\in\wedge^2\g$, denote by $r'$ and $r''$ the
left and right invariant bivector fields on $G$ corresponding to $r$.
Then it follows from \bref{cYB} that the bivector field $r'-r''$ defines a Poisson
bracket on $G$ which makes $G$ into a Poisson-Lie group.
On the other hand, a Poisson brackets on $M$ admitting,
in principle, a $\Uh$ invariant quantization endows $M$
with a structure of $(G,r)$-manifold. This means that
the action $G\times M\to M$ is a Poisson map.     
So, \Cor{coro2.3} describes the form of Poisson brackets on $M$
making $M$ into a $(G,r)$-manifold. One sees, in particular, that 
the classification of $(G,r)$ Poisson structures on $M$
reduces to the classification of $\ff$-brackets on $M$.

\spa We will also consider two parameter quantizations on $M$.
A two parameter quantization  of $\A$ is
an algebra $\A_{t,h}$ isomorphic to $\A[[t,h]]$
as a $\C[[t,h]]$ module and having a multiplication of the form
\be{}\label{twom}
m_{t,h}=m_0+t m_1^\prime+h m_1^{\prime\prime}+o(t,h).
\ee{}
With such a quantization one associates two Poisson brackets:
the bracket $v(a,b)=m^\prime_1(a,b)-m^\prime_1(b,a)$ along $t$,
and the bracket
$p(a,b)=m^{\prime\prime}_1(a,b)-m^{\prime\prime}_1(b,a)$ along $h$.
It is easy to check that $p$ and $v$ are compatible Poisson brackets,
i.e., their Schouten bracket $\[p,v\]$ is equal to zero.

\spa\begin{cor}{coro2.4} Let $\A_{t,h}$ be a $\Uh$ invariant
associative quantization of the form \bref{twom}.
Then the Poisson bracket  $p(a,b)=m_1^{\prime\prime}(a,b)-m_1^{\prime\prime}(b,a)$
has the form
\be{*}
p(a,b)=f(a,b)-r_M(a,b),
\ee{*}
where $r$ is the $r$-matrix
corresponding to $\Uh$ and $f$ is a $\ff$-bracket
with $\ff=\[r,r\]$. The Poisson bracket 
$v(a,b)=m_1^\prime(a,b)-m_1^\prime(b,a)$ in invariant
and compatible with $p$.
\end{cor}
\begin{proof} Similar to \Cor{coro2.3}.
\end{proof}

\spa In the following, a $\ff$-bracket on $M$ compatible with
a nondegenerate Poisson bracket we call a {\em good} bracket.

\section{Classification of $\ff$- and good brackets on semisimple orbits}

\pa 
Let $G$ be a complex connected simple Lie group with the Lie
algebra $\g$.
Let $\lf$ be a Levi subalgebra of $\g$, the Levi factor of a parabolic subalgebra.
Let $L$ be a Lie subgroup of $G$ with Lie algebra $\lf$. Such a subgroup is called
a Levi subgroup.
It is known that $L$ is a closed connected subgroup. Denote $M=G/L$ and
let $o\in M$ be the image of the unity by the natural projection $G\to M$.
Then $L$ is the stabilizer of $o$. It  is known, that $M$ may be
realized as a semisimple orbit of $G$ in the coadjoint representation $\g^*$. 
Conversely, any semisimple orbit in $\g^*$ is a quotient of $G$ by a Levi subgroup.  

\pa  
Let $\hh\subset \lf$ be a Cartan subalgebra of $\g$ and 
$\Omega_\lf\subset\Omega\subset \hh^*$ 
the sets of roots of $\lf$ and $\g$ corresponding to $\hh$.
Choose root vectors $E_\aa$, $\aa\in\Omega$, in such a way that
\be{}\label{ki}
(E_\aa,E_{-\aa})=1
\ee{} 
for the Killing form $(\cdot,\cdot)$ on $\g$.

\pa Let $Q$ be a set embedded in a linear space $V$ 
such that $0\not\in Q$ and $Q=-Q$.
We call a subset $B\subset Q$ {\em a linear subset} if
$B=Q\cap V_B$ where $V_B$ is the linear subspace in $V$
generated by $B$.
We call a subset $B\subset Q$ {\em semilinear}
if it follows from $x,y\in B$, $x+y\in Q$
that $x+y\in B$, and, in addition, $B\cap(-B)=\varnothing$, $B\cup(-B)=Q$.   
For a linear subset $B$ of Q,  we denote by $Q/B$ the
image of $Q$ without zero by the projection $F\to F/V_B$.

\pa Since $\lf$ is a Levi subalgebra, $\Omega_\lf$ is a linear subset
in $\Omega$. We put $\Ob=\Omega/\Omega_\lf$ and call
elements of $\Ob$ {\em quasiroots}.
For $\alpha\in\Omega$ we denote by $\alb$ its image in $\Ob$.
Let $Y$ be a semilinear subset in $\Ob$. One can easily shown
that there is a subset $P\subset Y$ such that any element of $Y$
can be uniquely presented as a linear combination of
elements of $P$ with integer coefficients.
We call $P$ a {\em  set of simple quasiroots} corresponding to $Y$
and $Y$ {\em a set of positive quasiroots} with respect to $P$.
It is clear that there is a set of simple roots, $\Pi$,
in $\Omega$ such that $P=\Pib$. Then $Y=\Ob^+$, where
$\Omega^+$ is the system of positive roots corresponding to $\Pi$. 
For such a $\Pi$ there is a subset, $\Ga\subset\Pi$, such
that $\lf$ coincides with the Lie subalgebra $\g_\Ga$, the subalgebra
generated by $\hh$ and elements $E_{\pm \aa}$, $\aa\in \Ga$.

\pa\label{ssec2.5} The projection $\pi:G\to M$ induces the map
$\pi_*:\g\to T_o$ where $T_o$ is the tangent space to $M$ at the
point $o$. Since the $\ad$-action of $\lf$ on $\g$ is semisimple,
there exists an $\ad(\lf)$-invariant subspace $\mm=\mm_\lf$ of $\g$
complementary to $\lf$, and one can identify $T_o$ and $\mm$ by means of
$\pi_*$. It is easy to see that subspace $\mm$ is uniquely defined
and has a basis consisting of the elements
$E_\gamma$, $\gamma\in\Omega\setminus\Omega_\lf$.

\pa\begin{propn}{prop2.1}The space $\mm$ considered as
a $\lf$ representation space decomposes into the direct sum of
subrepresentations $\mm_\beb$, $\beb\in \Ob$,
where $\mm_\beb$ is generated by all the elements $E_\beta$, $\beta\in \Omega$,
such that the projection of $\beta$ is equal to $\beb$. This decomposition
have the following properties:

a) all $\mm_\beb$ are irreducible;

b) $\mm_{-\beb}$ is dual to $\mm_\beb$;
 
c) for $\beb_1,\beb_2\in \Ob$ such
that $\beb_1+\beb_2\in \Ob$ one has
$[\mm_{\beb_1},\mm_{\beb_2}]=\mm_{\beb_1+\beb_2}$;

d) for any pair $\beb_1,\beb_2\in \Ob$  the representation 
$\mm_{\beb_1}\ot\mm_{\beb_2}$ is multiplicity free.
\end{propn}

\begin{proof}
Statements a), b), and c) are proven in \cite{DGS}, Remark 3.1. Statement d) follows
from the fact that the weight subspaces of all $\mm_\bbe$ have 
dimension one (see N.Bourbaki, Groupes et alg\`ebres de Lie,
Chap. 8, \S 9, Ex. 14).
\end{proof}

\pa Restricting to the point $o\in M$ defines the natural 
one-to-one correspondence between $G$ invariant
tensor fields on $M$ and $\lf$ invariant tensors over $\mm$.

Since $\lf$ contains a Cartan subalgebra, $\hh$, each $\lf$ invariant
tensor over $\mm$ is of weight zero with respect to $\hh$. It follows that there are
no invariant vectors in $\mm$. Hence, there are no invariant vector
fields on $M$.

\pa\begin{propn}{prop2.2}A bivector  $v\in\wedge^2\mm$ is $\lf$ invariant if and only if
it has the form
$v=\half\sum c(\alb)E_\alpha\wedge E_{-\alpha}$
where the sum runs over $\alpha\in \Omega\setminus\Omega_\lf$ 
(we suppose $c(-\baa)=-c(\baa)$).
\end{propn} 

\begin{proof} Follows from (\ref{ki}) and Proposition \ref{prop2.1}
(see also \cite{DGS}, Proposition 3.2).
\end{proof}

\pa Denote by $\[v,w\]\in\wedge^{k+l-1}\mm$ the Schouten bracket of
polyvector fields $v\in\wedge^k\mm$ and $w\in\wedge^l\mm$ defined
by the formula
\be{*}
\[X_1\wedge\cdots\wedge X_k, Y_1\wedge\cdots\wedge Y_l\]=\sum
(-1)^{i+j}[X_i,Y_j]_\mm\wedge X_1\wedge\cdots \hat X_i \cdots
\hat Y_j\cdots \wedge Y_l,
\ee{*}
where $[\cdot,\cdot]_\mm$ is the composition of the Lie bracket in $\g$ and
the projection $\g\to\mm$.
The defined Schouten bracket is compatible with the Schouten bracket on $M$
under identifying $\lf$ invariant polyvectors over $\mm$ and $G$ invariant
polyvector fields on $M$.

\pa It is obvious that any $\lf$ invariant bivector is $\theta$ 
anti-invariant for the Cartan automorphism $\theta$, $\theta(E_\alpha)=-\Ema$, of $\g$. 
Hence, if  $v,w\in\wedge^2\mm$ are $\lf$ invariant, then 
$\[v,w\]$ is $\theta$ invariant, i.e.,  is of the form 
$\[v,w\]=\sum e(\alpha,\beta)E_{\aa+\beta}\wedge E_{-\aa}\wedge E_{-\beta}$
where  $e(\aa,\beta)=-e(-\aa,-\beta)$.
Hence, in order to calculate  $\[v,w\]$ for such $v$ and $w$ 
it is sufficient to calculate
coefficients $e(\aa,\beta)$ for positive $\aa$ and $\beta$ by any choice
of the system of positive roots.  

\pa \begin{lemma}{lem2.3}
Let $v=\sum c(\aa)\Ea\wedge\Ema$, $w=\sum d(\aa)\Ea\wedge\Ema$
be elements from $\g\wedge\g$. Choose a system of positive roots.
Then for any positive roots $\aa,\beta,(\aa+\beta)$
the coefficient by the term 
$E_{\aa+\beta}\wedge E_{-\aa}\wedge E_{-\beta}$ in
$\[v,w\]$
is equal to
\be{}\label{coef}
N_{\aa,\beta}(d(\aa)(c(\beta)-c(\aa+\beta))+ 
d(\beta)(c(\aa)-c(\aa+\beta))-  
d(\aa+\beta)(c(\aa)+c(\beta))),   
\ee{}
where the number $N_{\aa,\beta}$ is defined by relation 
$[E_\aa,E_\beta]=N_{\aa,\beta}E_{\aa+\beta}$.
\end{lemma}

\begin{proof} Direct computation, see \cite{KRR}.
\end{proof}

\pa Let $\ff\in\wedge^3\g$ be an invariant element. Since $\g$ is simple,
$\ff$ is defined uniquely up to a factor. Denote by $\ff_M$ the invariant 
three-vector field on $M$ induced by $\ff$ with the help of the
action map $\g\to\Vect(M)$.
It is easy to check 
that $\ff_M$ is $\theta$ invariant and
up to a factor  has the form 
\be{}\label{ffM}
\ff_M=\third\sum_{\aa,\beta,\aa+\beta\in\Omega\setminus\Omega_\lf}
N_{\aa,\beta}E_{\aa+\beta}\wedge E_{-\aa}\wedge E_{-\beta}.
\ee{}

\pa From Lemma \ref{lem2.3} it follows that the Schouten bracket of
bivector $v=\half\sum c(\baa)E_{\aa}\wedge E_{-\aa}$ 
with itself is equal to $K^2\ff_M$ for a complex number $K$, if and only if
the following equations hold
\be{}\label{ff} 
c(\baa+\bbe)(c(\baa)+c(\bbe))=c(\baa)c(\bbe)+K^2
\ee{}
for all the pairs of quasiroots $\baa, \bbe$ such that 
$\baa+\bbe$ is a quasiroot.

Let  $X_{K^2}$ be the algebraic variety consisting of the points $\{c(\baa), \baa\in \Ob\}$
satisfying (\ref{ff}) (we always assume $c(\baa)=-c(-\baa)$).
So, for a given $\ff$ the variety $X=X_{K^2=-1}$ is the variety of all $\ff$-brackets.
It is clear that all the varieties $X_{K^2}$, $K\neq 0$, are isomorphic to $X$.

\pa Let $\{c(\baa)\}$ be a solution of (\ref{ff}) for a number $K$,
i.e.,  $\{c(\baa)\}\in X_{K^2}$.
It is easy to derive the following properties.

(*) If $c(\baa)+c(\bbe)=0$ then necessarily $c(\baa)=\pm K$, $c(\bbe)=\mp K$.

(**)   If $c(\baa)=\pm K$ and $c(\bbe)\neq \pm K$, then $c(\baa+\bbe)=\pm K$ and 
$c(\baa-\bbe)=\pm K$.

(***) If $c(\baa)=\pm K$ and $c(\bbe)=\pm K$, then $c(\baa+\bbe)=\pm K$.

\pa Formally, all the solutions of (\ref{ff}) for a fixed $K$ 
can be obtained in the following way. 
Choose a system of positive quasiroots, $\Ob^+$.
Denote by $\Pib$ the corresponding set of simple quasiroots.
Given $c(\baa)$ and $c(\bbe)$, we find from (\ref{ff}) that
\be{}\label{ab}
c(\baa+\bbe)=\frac{c(\baa)c(\bbe)+K^2}{c(\baa)+c(\bbe)}.
\ee{}
 Assume, $\baa,\bbe,\bga$ are positive quasiroots such that
$\baa+\bbe, \bbe+\bga, \baa+\bbe+\bga$ are also quasiroots.
Then the number $c(\baa+\bbe+\bga)$ can be calculated formally 
(ignoring possible division by zero) in two ways,
using (\ref{ab}) for  the pair $c(\baa), c(\bbe+\bga)$ on the right hand side
and  also for the pair $c(\baa+\bbe), c(\bga)$.
But it is easy to check that these two ways
give the same value of $c(\baa+\bbe+\bga)$. In this sense the system
of equations corresponding to  (\ref{ab}) for all pairs is consistent.
So, taking arbitrary values $c(\baa)$ for simple quasiroots $\baa$ one can 
try to find $c(\baa)$ for all $\baa\in\Ob^+$ recursively. 
We say that a solution, $\{c(\baa)\}$, of \bref{ff} can be obtained recursively 
if in the course of the recursive procedure started with the values
$c(\baa)$ for simple quasiroots $\baa$ the denominators in (\ref{ab})
will be not equal to zero. 

\pa \begin{propn}{prop2.5} For $K\neq 0$ the following holds.

a) Any solution of (\ref{ff}) can be obtained recursively
by choosing a respective system of positive quasiroots.

b) The variety $X_{K^2}$ is without singularities, connected, and of dimension $k$,
where $k$ is equal to the number of simple quasiroots.
\end{propn}

\begin{proof}
For proving a) we have to show that for any solution
 $\{c(\baa)\}\in X_{K^2}$ one can choose a system of positive
quasiroots in such a way that the denominators appearing in (\ref{ab})
by the recursive procedure are
not equal to zero. It follows  from (**) that the set $\Psi$ consisting of $\baa$ such that
$c(\baa)\neq\pm K$ is a linear subset of $\Ob$. Moreover, the function
$c(\baa)$ is constant on the cosets of $\Ob/\Psi$.
Let $Y$ is the set of cosets on which this function has the value $K$.
It follows from (***) that $Y$ is a semilinear subset of $\Ob/\Psi$.
Let $\Ob^+$ be a semilinear subset of $\Ob$ projecting on $Y$.  
Then it follows from (*) that for $\baa,\bbe\in\Ob^+$
$c(\baa)+c(\bbe)\neq 0$, which proves a).

Let $\Pib=\{\baa_i, i=1,...,k\}$ be the set of simple quasiroots corresponding to $\Ob^+$.
Let $c(\baa_i)=c_i$. It is clear that 
starting the recursive procedure with 
$c(\baa_i)=c_i^\prime$ for $c_i^\prime$ arbitrary but
close enough to $c_i$, the denominators in (\ref{ab})
remain not equal to zero. This proves that any point of $X_{K^2}$ 
is non-singular and $X_{K^2}$ has dimension $k$.

Let us prove that $X_{K^2}$ is connected. Fix a set of positive quasiroots, $\Ob^+$,
and the corresponding set of simple quasiroots, $\Pib=\{\baa_i, i=1,...,k\}$. 
We say that a $k$-tuple of complex numbers $(c_1,...,c_k)$ is admissible, if
starting with $c(\baa_i)=c_i$ one obtains a solution of (\ref{ff})
by the recursive procedure. 
It is clear that the admissible tuples form 
a subset, $A$, of $\C^k$ complement
to an algebraic subset of lesser dimension, therefore $A$ is connected.
On the other hand, the set of points $\{c(\baa)\}\in X_{K^2}$ such that
$c(\baa_i)$ form an admissible $k$-tuple is obviously dense in $X_{K^2}$.
This proves the connectness of $X_{K^2}$.
\end{proof}

\pa 
Let $\{c(\baa)\}$ be a solution of (\ref{ff}) and 
$\Psi\in\Ob$ the linear subset of $\Ob$ such that 
$c(\baa)\neq\pm K$ for $\baa\in\Psi$. Then, using the formula
for $\coth(x+y)$, similar to (\ref{ab}), one can see
that there exists a linear form $\lambda:\Psi\to\C$ such that
\be{}\label{abla}
\lambda(\baa)\not\in\frac{2\pi i}{K}\ZZ \qquad \baa\in\Psi
\ee{}
and
\be{}\label{ab1}
c(\baa)=K\coth\(\frac{K}{2}\lambda(\baa)\) \qquad \mbox{for any} \quad\baa\in\Psi.
\ee{}  

As $K\to 0$, (\ref{abla}) makes into
\be{}\label{ablalim}
\lambda(\baa)\neq 0 \qquad \baa\in\Psi
\ee{}
and (\ref{ab1}) tends to 
\be{}\label{ab1lim}
c(\baa)=\frac{1}{\lambda(\baa)} \qquad \mbox{for any} \quad\baa\in\Psi.
\ee{}
So, we come to 

\pa\begin{propn}{prop2.6}
a) For $K\neq 0$, any solution of (\ref{ff}) is determined by: 
choosing a linear subset, $\Psi$, in $\Ob$, a semilinear subset, $B$, in $\Ob/\Psi$,
and a linear form, $\lambda:\Psi\to\C$, satisfying (\ref{abla}). 
The respective solution $\{c(\baa)\}$ is the following: 
for $\baa\in \Psi$ $c(\baa)$ is defined by (\ref{ab1});
for $\baa\not\in \Psi$, $c(\baa)= K$ if
the projection of $\baa$ in $\Ob/\Psi$ belongs to $B$,
$c(\baa)=-K$ if
the projection belongs to $-B$.

b) For $K=0$, any solution of (\ref{ff}) is defined by: 
choosing a linear subset, $\Psi$, in $\Ob$
and a linear form, $\lambda:\Psi\to\C$, satisfying (\ref{ablalim}). 
The solution $\{c(\baa)\}$ is the following: 
for $\baa\in \Psi$ $c(\baa)$ is defined by (\ref{ab1lim});
for $\baa\not\in \Psi$ $c(\baa)= 0$.
\end{propn}

\pat{\em Remark.} As mentioned in \cite{Lu}, for $\lf=\hh$ the solutions
described in \Prop{prop2.6} 
relate to solutions of classical dynamical Yang-Baxter equations, \cite{EV}.

\pa Note that by $K=0$ the solutions of (\ref{ff}) define
Poisson brackets on $M$, so \Prop{prop2.6} b) describes all the Poisson
brackets on $M$. We see that nondegenerate Poisson brackets on $M$
are in one-to-one correspondence with  the linear forms
$\lambda:\Ob\to \C$ such that $\lambda(\baa)\neq 0$ for all $\baa\in\Ob$
and have the form
\be{}\label{kirbr}
\half\sum_{\baa\in\Ob}\frac{1}{\lambda(\baa)}E_\aa\wedge \Ema.
\ee{} 
This is exactly the KKS bracket on the orbit in $\g^*$ passing through
the linear form on $\g$ being the trivial extension of $\lambda$. 

\pa\label{ssec2.21} Denote by $X_0$ the variety of nondegenerate Poisson brackets on $M$.
Since $X_0$ coincides with all solutions of \bref{ff}, $\{c(\baa)\}$, 
such that $\Pi_\baa c(\baa)\neq 0$, it is clear that $X_0$ is
an affine connected algebraic variety without singularities.

\pa Fix a Poisson bracket $s$ of the form \bref{kirbr}.
Let us describe the invariant brackets 
$f=\sum c(\baa)E_{\aa}\wedge E_{-\aa}$ satisfying 
the conditions
\be{}\label{pbf}
&\[f,f\]=K^2\ff_M, \\ 
&\[f,s\]=0. \label{pbfc} 
\ee{}
with some $K\neq 0$.

A direct computation shows that  conditions \bref{pbf} and \bref{pbfc} are equivalent 
to the system of equations  for the coefficients $c(\baa)$ of $f$, \cite{DGS},
\be{} \label{eq1}
c(\bbe)\ll(\bbe)=c(\baa)\ll(\baa)\pm K\ll(\baa+\bbe)\\
c(\baa+\bbe)\ll(\baa+\bbe)=c(\baa)\ll(\baa)\pm K\ll(\bbe) \label{eq2}
\ee{}
with the same sign before $K$, for all the pairs of quasiroots $\baa, \bbe$ such that 
$\baa+\bbe$ is a quasiroot.

\pat{Definition.}
Let $M$ be an orbit in $\g^*$ (not necessarily semisimple).
The invariant bracket $f$ on $M$ is said to be {\em good}
if $f$ satisfies conditions (\ref{pbf}) and (\ref{pbfc}) for $s$ 
the Kirillov-Kostant-Souriau (KKS)
Poisson bracket on $M$.
We call $M$ a {\em good} orbit, if there exists a good bracket on it.

\pa\begin{propn}{prop2.7}

a) For $\g$ of type $A_n$ all semisimple orbits are good.

b) For  all other $\g$, the orbit $M$ is good if and only if $\lf=\g_\Ga$
where $\Ga\subset\Pi$ for a system $\Pi$ of simple root for $\g$ and
the set $\Pi\setminus\Ga$ consists of one or two roots which appear   
in the representation of the maximal root with coefficient 1.

c) For a given $K\neq 0$ the good brackets $f$ on a good orbit form a one-dimensional variety:
all such brackets have the form 
\be{*}
\pm f_0+ts,
\ee{*}
where $t\in \C$ and $f_0$ is a fixed bracket satisfying 
(\ref{pbf}), (\ref{pbfc}).
\end{propn}

\begin{proof} The proof reduces to solving the system of
equations defined by \bref{pbf} and \bref{pbfc}, see \cite{DGS}.
\end{proof}

\pa So, if the set $\Pi\setminus\Ga$ consists of one root, $M$ is exactly
a hermitian symmetric space. As follows from the classification
of simple Lie algebras, the case when the set $\Pi\setminus\Ga$ consists of two
roots appears (besides $A_n$) for $\g$ of types $D_n$ and $E_6$.

\pa\label{ss3.26}
\Prop{prop2.7} shows that the property for $M$ to be a good orbit
depends only on the pair $(\g,\lf)$ but not on the realization of $M$ as
an orbits. The pair $(\g,\lf)$ corresponding to a good orbit we
call a {\em good} pair.

\pat{\em Remark.}\label{rem2.1} It is clear that if $f$ satisfies \bref{pbf} and \bref{pbfc} then
$\pm f+ts$ also satisfies the same conditions (with the same $K$) for all numbers $t$.
\Prop{prop2.7} c) shows that, conversely, all good brackets on a good semisimple
orbit are contained in these families
$\pm f+ts$, $t\in\C$. 

Denote by $Y$ the variety of good brackets $f$ on $M$
satisfying \bref{pbf} and \bref{pbfc} for a fixed $K\neq 0$ and
some nondegenerate Poisson bracket $s$.
From the above it follows that there is a projection, $Y\to X_0$,
$f\mapsto s$, where $s$ is a Poisson bracket such that $\[f,s\]=0$.
The fiber over $s\in X_0$ consists of two components,
$\{\pm f_0+ts\}$, $t\in\C$, isomorphic to $\C$.
These components correspond to choosing the sign in
\bref{eq1}, \bref{eq2}.

\pat{\em Remark.} It is shown in \cite{Do1} that in case $A_n$, i.e., when
$\g=sl(n)$, all the coadjoint orbits (not necessarily semisimple) are good.
Moreover, there exists a unique quadratic $\ff$-bracket
on $sl(n)^*$ which can be restricted to all orbits to give
good brackets on them. This quadratic bracket on $sl(n)^*$ can be quantized,
\cite{Do2}.

\section{Poisson complexes}
\pa Let $C^k=(\Lambda^k\mm)^\lf$ be the space of $\lf$ invariant $k$-vectors
on $\mm$. This space is identified with the space of $G$ invariant
$k$-vector fields on $M$.
Denote by $\CC^k$ the sheaf of holomorphic functions on $X_{K^2}$ with values in $C^k$.
We form the complex $(C^\bullet, \delta_f)$ where $\delta_f$
is the differential given by the Schouten
bracket with a bivector $f\in X_{K^2}$,
$$\delta_f:u\mapsto \[f,u\]\,\qquad\mbox{for} 
\quad u\in C^\bullet. $$
The condition  $\delta_f^2=0$ follows from the Jacobi identity for
the Schouten bracket together with the fact that $\[f_x,f_x\]=K^2\ff_M$.

We also consider the complex of sheaves $(\CC^\bullet,\delta)$ on
$X_{K^2}$. The operator $\delta$ is defined as  
\be{}\label{cob}
\delta(u)(f)=\[f,u(f)\]=\delta_f(u(f)),
\ee{}
where $u$ is a section of $\CC^\bullet$ and $f\in X_{K^2}$.

We denote by $H^k(M,\delta_f)$ and $H^k(\CC^\bullet,\delta)$ the
cohomologies of $(C^\bullet, \delta_f)$ and $(\CC^\bullet,\delta)$ , 
respectively, whereas the usual
de Rham cohomologies are denoted by $H^k(M)$.

\pa\begin{propn}{prop3.1}
a) For any nondegenerate Poisson brackets $f\in X_0$ 
and, if $K\neq 0$,
for almost all $f\in X_{K^2}$ (except an algebraic subset
of lesser dimension) one has
\be{}\label{a}
H^k(M,\delta_f)=H^k(M)
\ee{}
for all $k$.
In particular, $H^k(M,\delta_f)=0$ for odd $k$. 

b) Let $K\neq 0$. Then $H^2(M,\delta_f)=H^2(M)$ for all $f\in X_{K^2}$.
\end{propn}

\begin{proof} The proof of a) follows \cite{DGS}.
First, let $v$ be a nondegenerate Poisson bracket on $M$, in particular, $v\in X_0$.
Then the complex of polyvector fields on $M$, $\Theta^\bullet$, with
the differential $\delta_v$ is well defined.
Denote by $\Omega^\bullet$ the de Rham complex on $M$.
Since none of the coefficients $c(\baa)$ of $v$ are
zero, $v$ is a nondegenerate bivector field, and therefore
it defines an $\A$-linear isomorphism 
$\tilde{v}:\Omega^1\to\Theta^1$, $\omega\mapsto v(\omega,\cdot)$,
which can be extended up to the isomorphism $\tilde{v}:\Omega^k\to\Theta^k$
of $k$-forms onto $k$-vector fields for all $k$.
Using Jacobi identity for $v$ and invariance of $v$, one can
show that $\tilde{v}$ gives a $G$ invariant isomorphism of
these complexes, so their cohomologies are the same.

Since $\g$ is simple, the subcomplex of 
$\g$ invariants, $(\Omega^\bullet)^\g$, splits off
as a subcomplex of $\Omega^\bullet$.
In addition, $\g$ acts trivially on cohomologies,
since for any $g\in G$ the map $M \to M$,
$x\mapsto gx$, is homotopic to the identity map, ($G$ is a connected
Lie group corresponding to $\g$).
It follows that cohomologies of complexes
$(\Omega^\bullet)^\g$ and $\Omega^\bullet$ coincide.

But $\tilde{v}$ gives an isomorphism of complexes
$(\Omega^\bullet)^\lf$ and 
$(\Theta^\bullet)^\g=((\Lambda^\bullet\mm)^\lf,\delta_v)$.
So, cohomologies of the latter complex coincide with de Rham
cohomologies, which
proves a) for $v$ being Poisson brackets.

Denote by $X$ the variety of points $\{c(\baa),\baa\in\Ob;K\}$
satisfying \bref{ff}. Consider the family of complexes
$((\Lambda^\bullet\mm)^\lf,\delta_v)$, $v\in X$. It is clear that
$\delta_v$ depends algebraicly on $v$. It follows from the uppersemicontinuity  
of $\dim H^k(M,\delta_v)$ and the fact that $H^k(M)=0$ for odd $k$,
\cite{Bo}, that $H^k(M,\delta_v)=0$ for odd $k$ and almost all $v\in X$.
Using the uppersemicontinuity again and the fact that
the number $\sum_k(-1)^k\dim H^k(M,\delta_v)$ is the same for all $v\in X$,
we conclude that $\dim H^k(M,\delta_v)=\dim H^k(M)$ for even $k$ and
almost all $v\in X$. We have that there is a number $K_0\neq 0$ and
$v\in X_{K_0}$ such that \bref{a} holds for $\delta_v$. Since the cohomologies
do not change when $v$ replaces by $cv$ with any number $c\neq 0$,
we conclude that for any $K\neq 0$ there is $f$ satisfying \bref{a}.
Since by \Prop{prop2.5} b) $X_{K^2}$ is connected, it follows that
\bref{a} holds for almost all $f\in X_{K^2}$.

Let us prove b). Since $(\Lambda^1\mm)^\lf=0$, $H^2(M,\delta_f)$ coincides
with $\Ker(\delta_f:C^2\to C^3)$.
According to \Prop{prop2.5},
let us choose a system of positive quasiroots $\Ob^+$ such that $f$ 
has the form $f=\sum_{\baa\in\Ob^+} c(\baa) E_{\alpha}\wedge \Ema$, 
and for $\baa,\beb,\baa+\beb\in\Ob^+$ $c(\baa)+c(\beb)\neq 0$. 
By \Lem{lem2.3}, any $v\in C^2$ such that $\[f,v\]=0$ has the form
$v=\sum_{\baa\in\Ob^+} d(\baa) E_{\alpha}\wedge \Ema$
where $d(\baa)$ obey the equation
\be{}\label{defd}
d(\aa)(c(\beta)-c(\aa+\beta))+ 
d(\beta)(c(\aa)-c(\aa+\beta))-  
d(\aa+\beta)(c(\aa)+c(\beta))=0.
\ee{}
Starting with arbitrary values $d(\baa)$ for simple quasiroots we find
$d(\baa)$ for all $\baa\in\Ob^+$ recursively using the formula
following from \bref{defd}
\be{}\label{defd1}
d(\baa+\beb)=\frac{d(\aa)(c(\beta)-c(\aa+\beta))+ 
d(\beta)(c(\aa)-c(\aa+\beta))}{c(\aa)+c(\beta)}.
\ee{}
Indeed, denominators of \bref{defd1} are not equal to zero by choosing
of $\Ob^+$. Moreover,
assume $\baa,\bbe,\bga$ are positive quasiroots such that
$\baa+\bbe, \bbe+\bga, \baa+\bbe+\bga$ are also quasiroots.
Then the number $d(\baa+\bbe+\bga)$ can be calculated  in two ways,
using (\ref{defd1}) for  the pair $d(\baa), d(\bbe+\bga)$ on the right hand side
and  also for the pair $d(\baa+\bbe), d(\bga)$.
But it is easy to check that these two ways
give the same value of $d(\baa+\bbe+\bga)$.

So we see that 
$\dim H^2(M,\delta_f)=$ \{{\em the number of simple quasiroots}\}
for all $f\in X_{K^2}$.
\end{proof}

\pa\begin{propn}{prop3.2}Let $K\neq 0$. Then
\be{*}
H^3(\Ga\CC^\bullet,\delta)=0,
\ee{*}
where $\Ga\CC^k$ is the space of global sections of $\CC^k$ over $X_{K^2}$.
\end{propn}

\begin{proof}
Since $X_{K^2}$ is a Stein manifold, it is enough to prove
that $H^3(\CC^\bullet,\delta)=0$.
Note that according to \Prop{prop3.1} b)
$\Ker(\delta:\CC^2\to\CC^3)=H^2(\CC^\bullet,\delta)$ is
a subbundle (direct subsheaf) of $\CC^2$. Therefore, $\Im(\delta)$ is
a subbundle of $\CC^3$. On the other hand, 
$\Ker(\delta:\CC^3\to\CC^4)$ being
a subsheaf of $\CC^3$ is a torsion free sheaf. It follows that 
$H^3(\CC^\bullet,\delta)$ is a torsion free sheaf. According to
\Prop{prop3.1} a) the support of $H^3(\CC^\bullet,\delta)$ is
an algebraic subset of $X_{K^2}$ of lesser dimention. Hence
$H^3(\Ga\CC^\bullet,\delta)=0$.
\end{proof}

\pa\label{sec3.3} Let $M$ be a good orbit, $s$ the KKS Poisson bracket on $M$, and $f_0$
a good bracket on $M$. Hence, $\[f_0,f_0\]=K^2\ff_M$ for some $K\neq 0$
and $\[f_0,s\]=0$. Consider the family of brackets
$$\phi=f_{h,t}=hf_0+ts, \qquad h,t\in\C.$$ 
One has $\[f_{h,t},f_{h,t}\]=h^2K^2\ff_M$,
therefore, according to \Prop{prop2.7} c) this family
contains all the good brackets on $M$ for all $K$.
We will consider $f_{h,t}$ as a linear map of $\C^2$
with the coordinates $(h,t)$ to the space $(\Lambda^2\mm)^\lf$.

Denote by $\CC^k(m)$ the space of homogeneous maps 
$\C^2\to(\Lambda^k\mm)^\lf$ of degree $m$. Thus,
$f_{h,t}\in \CC^2(1)$. Define the differential 
$\delta:\CC^\bullet(m)\to\CC^\bullet(m+1)$ as
$(\delta(a))(h,t)=\[f_{h,t},a(h,t)\]$.

\pa\begin{propn}{prop3.3}a) For all $(h,t)\in\C^2\setminus 0$ 
$$H^2(M,\delta_{f_{h,t}})=H^2(M).$$ 

b) Let $m\geq 0$ and $b\in\CC^3(m+1)$ such that $\delta b=0$.
Then there exists $a\in\CC^2(m)$ such that $\delta a=b$.
\end{propn}

\begin{proof}Let us prove a).
If $h\neq 0$ then $f_{h,t}\in X_{hK}$ with $hK\neq 0$. Hence
$H^2(M,\delta_{f_{h,t}})=H^2(M)$ by \Prop{prop3.1} b).
If $h=0$ then $t\neq 0$ and $f_{0,t}$ is a nondegenerate
Poisson bracket on $M$, therefore,
$H^2(M,\delta_{f_{0,t}})=H^2(M)$ by \Prop{prop3.1} a). 

Let us prove b). Denote by $\L(m)$ a linear holomorphic vector bundle of
degree $m$ over the Riemann sphere $\SS$. The space of global
sections of $\L(m)$ may be naturally identified with the space
of homogeneous polynomials of two variables of degree $m$.
Taking as the variables $h$ and $t$, one can consider the space $\CC^k(m)$ as the space
of global sections of a vector bundle, $\EE^k(m)$, that is a direct sum
of $\dim(\Lambda^k\mm)^\lf$ copies of $\L(m)$. 
It is obvious that $f_{h,t}$ defines a map of sheaves,
$\delta^k_m:\EE^k(m)\to\EE^{k+1}(m+1)$. Denote
$\HH^k(m)=\Ker(\delta^k_m)/\Im(\delta^{k-1}_{m-1})$.

It follows from  a) that 
$\delta^2_m:\EE^2(m)\to\EE^{3}(m+1)$ is a
map of bundles, i.e., its image is a direct subsheaf of $\EE^3(m+1)$.
Hence, $\HH^3(m+1)$ is
a sheaf over $\SS$ without torsion. But over a neighborhood of the point of $\SS$
with homogeneous coordinates $(0,t)$  $\HH^3(m+1)=0$. This follows from
the fact that $f_{0,t}$ is a nondegenerate Poisson bracket
and, therefore, $H^3(M,\delta_{f_{0,t}})=0$. Since $\SS$ is connected, 
$\HH^3(m+1)=0$ over $\SS$. So, $\Ker\delta^3_{m+1}=\Im\delta^2_m$.

Consider the exact sequence of sheaves over $\SS$:
\be{*}
0\tor{}\HH^2(m)\tor{}\EE^2(m)\tor{}\Im{\delta^2_m}\tor{}0.
\ee{*}
To complete the proof of b) one needs to show that any
global section of $\Im{\delta^2_m}$ can be lifted to a global
section of $\EE^2(m)$. But this follows from the fact
that $H^1(\SS,\HH^2(m))=0$. To prove the last fact we observe that
$\HH^2(m)$ is a subbundle of $\EE^2(m)$, therefore,
is a direct sum of a number of copies of $\L(m)$.
Since $m\geq 0$, $H^1(\SS,\L(m))=0$. This implies that
$H^1(\SS,\HH^2(m))=0$, too.   
\end{proof}

\section{The $G$ invariant $\Phi_h$ associative quantization
in one parameter}

\pa Denote by $(HC^\bullet_M,\partial)$ the Hochschild complex on $M$,
where each space $HC^k_M$ consists of holomorphic
$k$-differential operators on $M$.
Let $X$ be a complex analytic manifold. The map $\psi:X\to HC_M^k$ is
called to be holomorphic if for any open subsets $U\subset X$ and
$V\subset M$ and any holomorphic functions
$a(x,y)_1,...,a(x,y)_k$ on $U\times V$
$\psi(a_1,...,a_k)$ is also a holomorphic function on $U\times V$.
We will denote the map $\psi$ by $\psi_x$, $x\in X$, and call
$\psi_x$ a holomorphic family of $k$-differential operators on $M$. 

Denote by $HC^k_M(X)$ the space of all holomorphic maps
$X\to HC_M^k$ and by $(HC^\bullet_M(X),\partial)$ the corresponding
complex. It is clear that $(HC^\bullet_M(X),\partial)$ is naturally identified
with the subcomplex of $(HC^\bullet_{X\times M},\partial)$ 
consisting of polydifferential operators along $M$. 
Denote by $\Lambda^kT_M(X)$ the space of analytic maps of $X$
to the space of polyvector fields on $M$.

\pa\begin{propn}{prop4.1}Let $X$ be a Stein manifold.
Then
\be{*}
H^k(HC^\bullet_M(X),\partial)=\Lambda^kT_M(X).
\ee{*}
\end{propn}

\begin{proof} The proof can be proceed in the similar way as
for $(HC^\bullet_M,\partial)$, using that $M$ and $X$ are Stein
manifolds, \cite{HKR}.
\end{proof}

\pa\begin{propn}{prop4.2}Let $\g$ be a simple Lie algebra,
$M$ a semisimple orbit in $\g^*$.
Let  $X=X_{K^2=-1}$ be the manifold of invariant
brackets $f$ satisfying $\[f,f\]=-\ff_M$.
Then there exists a holomorphic family of multiplications $\mu_{f,h}$ on $\A$
of the form
\be{}\label{aspr}
\mu_{f,h}(a,b)=ab +(h/2) f(a,b) +\sum_{n\geq 2}h^n \mu_{f,n}(a,b), \quad f\in X,
\ee{}
that  is $\Ug$ invariant  and
$\Phi_h$ associative.

\end{propn}
\begin{proof} The proof is essentially follows to \cite{DGS}, Proposition 5.1,
but here we construct the multiplication for all $f$ simultaneously 
using parameterized Poisson cohomologies from the previous section.

To begin,  consider the multiplication $\mu^{(1)}(a,b)=ab +(h/2) f(a,b)$.
The corresponding  obstruction cocycle  is given by 
$$obs_2=\frac1{h^2}(\mu^{(1)}(\mu^{(1)}\otimes id)-
\mu^{(1)}(id\otimes\mu^{(1)})\Phi_h)$$
considered modulo terms of order $h$. No $\frac1h$ terms appear because
$f$ is a biderivation and, therefore, a Hochschild cocycle. 
The fact that the presence of
$\Phi_h$ does not interfere with the cocyle condition and 
that this equation defines a Hochschild $3$-cocycle 
was proven in \cite{DS1} (see the proof of Proposition 4.1 there). 
By \Prop{prop4.1} the
differential Hochschild cohomology of $\A$ in dimension $p$ is 
the space of holomorphic families of $p$-polyvector fields on $M$
parameterized by $X$. 
Since $\g$ is reductive, the subspace of $\g$ invariants splits off
as a subcomplex and has cohomology given by $(\Lambda^p\mm)^{\lf}(X)$.
The complete antisymmetrization of a $p$-tensor projects
the  space of invariant differential $p$-cocycles onto the subspace
$(\Lambda^p\mm)^\lf(X)$ representing the cohomology. The equation
$\[f,f\]+\ff_M=0$ implies that the obstruction cocycle is a coboundary,
and we can find a $2$-cochain $\mu_{f,2}$, so that 
$\mu^{(2)}=\mu^{(1)}+ h^2\mu_{f,2}$ satisfies 
$$\mu^{(2)}(\mu^{(2)}\otimes id)-
\mu^{(2)}(id\otimes\mu^{(2)})\Phi_h=0\mbox{  mod }h^3.$$
Assume we have defined the deformation $\mu^{(n)}$ to order $h^n$ 
such that $\Phi_h$ associativity holds modulo $h^{n+1}$, then we define the
$(n+1)$-st obstruction cocycle by
$$obs_{n+1}=\frac1{h^{n+1}}(\mu^{(n)}(\mu^{(n)}\otimes id)-
\mu^{(n)}(id\otimes\mu^{(n)})\Phi_h)\mbox{ mod }h.$$

In \cite{DS1} (Proposition 4.1) it is shown that the usual 
proof that the obstruction cochain satisfies the cocycle condition
carries through to the $\Phi_h$ associative case. 
The coboundary of $obs_{n+1}$ appears as the 
$h^{n+1}$ coefficient of the signed sum of the compositions of 
$\mu^{(n+1)}$ with $obs_{n+1}$. 
The fact that  $\Phi_h=1$ mod $h^2$ together with  the pentagon
identity implies that the sum vanishes
identically, and thus all coefficients vanish, including the coboundary
in question. Let $obs_{n+1}'\in (\Lambda^3\mm)^\lf(X)$ be 
the projection of $obs_{n+1}$ on the totally skew symmetric part, which represents
the cohomology class of the obstruction cocycle. The coefficient of
$h^{n+2}$ in the same signed sum, when projected on the skew symmetric
part, is $\[f,obs_{n+1}'\]$ which is the coboundary of $obs'_{n+1}$  
in the complex  $((\Lambda^\bullet\mm)^\lf(X),\delta= \[f,.\])$. Thus
$obs'_{n+1}$  is a $\delta$ cocycle.  By Proposition \ref{prop3.2},
this complex has zero 3-cohomology. Now we  modify 
$\mu^{(n+1)}$ by adding a term $h^n\mu_{f,n}$ with $\mu_{f,n}\in
(\Lambda^2\mm)^\lf$ and consider the $(n+1)$-st obstruction 
cocycle for $\mu^{\prime(n+1)}=\mu^{(n+1)}+h^n\mu_{f,n}$. Since the term we added at
degree $h^n$ is a Hochschild cocyle, we do not introduce a $h^n$ term
in the calculation of $\mu^{(n)}(\mu^{(n)}\otimes id)-
\mu^{(n)}(id\otimes\mu^{(n)})\Phi_h $ and the totally skew symmetric
projection $h^{n+1}$ term has been  modified
by $\[f,\mu_{f,n}\]$. By choosing $\mu_{f,n}$ appropriately,
we can make the $(n+1)$-st obstruction cocycle represent the
 zero cohomology class, and we are able to continue
the recursive construction of the desired deformation.
\end{proof}

\pa Let $X$ be as in \Prop{prop4.2}. Let $(o,\C[[h]])$ be the formal manifold
which is the formal neighborhood of $0\in\C$. We call a morphism 
$\pi: o\to X$ a {\em formal path} in $X$. A formal path $\pi$ may be
given by a formal series in $h$, 
$$f(h)=f_0+hf_1+h^2f_2+\cdots, \qquad f_k\in(\Lambda^2\mm)^\lf$$
satisfying $\[f(h),f(h)\]=-\ff_M$. It is clear that $f_0\in X$;
we call it the {\em origin} of the path $f(h)$. The element $f_1$
belongs to the tangent space to $X$ at the point $f_0$, which
consists of elements $v$ such that $\[f_0,v\]=0$.

It is clear that if we put
$f=f(h)$ in the multiplication $\mu_{f,h}$, 
we obtain a $\Ug$ invariant $\Phi_h$ associative multiplication 
which depends only on $h$
and has the form
$$\mu_{f(h),h}(a,b)=ab+(h/2)f_0(a,b)+\cdots$$  
with $\ff_M$-bracket $f_0$. 
In particular, we obtain 

\pa\begin{cor}{cor4.2}Any $\ff_M$-bracket on a
semisimple orbit in $\g^*$ can be quantized.
\end{cor}

\begin{proof} Indeed, let $f_0$ be a $\ff_M$-bracket.
The multiplication corresponding to any path
$f(h)=f_0+hf_1+\cdots$ with origin in $f_0$ is as required.
\end{proof}

\pa\begin{propn}{prop4.2a}The multiplication $\mu_{f,h}$
of the form \bref{aspr}  has the following universal properties:

a) For any  $\Ug$ invariant $\Phi_h$ associative multiplication
$m_h$, there exists a formal path in $X$, $f(h)$, such that
$m_h$ is equivalent to $\mu_{f(h),h}$.

b) Multiplications corresponding to different paths are
not equivalent.
\end{propn}

We will need the following

\spa\begin{lemma}{lem4.2}Let
$m_h(a,b)=ab+hf_0(a,b)+h^2m_2(a,b)+\cdots$ 
be a multiplication and $f(h)$ a path such that
the multiplication
$\mu_{f(h),h}$ coincides with $m_h$ modulo $h^{n+1}$.
Then there exist a path $f'(h)=f(h)+h^{n}p_1+\cdots$
and a differential operator $D$ such that
the multiplication $m'_h=(1+h^{n+1}D)\circ m_h$,
$$m'_h(a,b)=(1+h^{n+1}D)^{-1}m_h((1+h^{n+1}D)a,(1+h^{n+1}D)b),$$
coincides with $\mu_{f'(h),h}$ modulo $h^{n+2}$.
\end{lemma}

\begin{proof}
We have 
$$\mu_{f(h),h}=ab+hf_0+h^2m_2+\cdots+h^nm_n+h^{n+1}\mu_{n+1}+\cdots,$$
where $m_k$, $k=2,3,...$, are terms appearing in the expansion of $m_h$.
It is easy to check that $\mu_{n+1}-m_{n+1}$ is a Hochschild cocycle,
because both $\mu_{n+1}$ and $m_{n+1}$ resolve the same obstruction
$obs_{n+1}$,
\be{*}
obs_{n+1}(a,b,c)=
\sum_{\stackrel{i+j=n}{ i,j\geq 1}}\(m_i(m_j(a,b),c)-m_i(a,m_j(b,c))\),
\ee{*}
 depending only on $m_k$, $k\leq n$.
Hence, one has $\mu_{n+1}=m_{n+1}+\partial D+p_1$, where
$D$ is a Hochschild 1-chain, i.e., a differential operator,
and $p_1$ is a bivector field. 
Applying $1+h^{n+1}D$ to $m_h$ we obtain
$$(1+h^{n+1}D)\circ m_h=\mu_{f(h),h}-h^{n+1}p_1 \quad \mbox{mod } h^{n+2}.$$
Observe now that $p_1$ is a $\delta_{f_0}$
cocycle, i.e., $\[f_0,p_1\]=0$. This follows from the fact
that $(1+h^{n+1}D)\circ m_h$ is a $\Phi_h$ associative multiplication.
Indeed, if $\[f_0,p_1\]$ is not equal to zero, its
contribution to $obs_{n+2}$ is not a Hochschild
coboundary.
So, $p_1$ is a tangent vector to $X$ at $f_0$.
Since, by \Prop{prop2.5} b), $X$ is without singularities,
there exists a formal path in $X$ of the form
$p(h)=f_0+h^{n}p_1+\cdots$. Let $f'(h)$ be the path
that in local coordinates on a neighborhood of $f_0$ in $X$ is the sum $f(h)+p(h)$.
It is clear that $f'(h)=f(h)+h^{n}p_1+\cdots$. Putting in
$\mu_{f,h}$  $f'(h)$ instead $f(h)$ does not change 
the coefficients by $h^k$, $k\leq n$, in $\mu_{f(h),h}$ and changes the $(n+1)$-st
coefficient adding $p_1$ to it. So we have
$$m'_h=(1+h^{n+1}D)\circ m_h=\mu_{f'(h),h} \quad \mbox{mod } h^{n+2},$$
as required 
\end{proof}

\spat{Proof of \Prop{prop4.2a}.} Let us prove a).
Let $f_0$ be the $\ff_M$-bracket
corresponding to the multiplication $m_h$. By \Cor{cor4.2}
one can assume that $m_h$ coincides modulo $h^2$
with $\mu_{f(h),h}$ for the trivial path $f(h)=f_0$.
Using \Lem{lem4.2} we find $D_2$ and $p_1$ such that the multiplication
$m^{(2)}_h=(1+h^2D_2)\circ m_h$ is equal modulo $h^3$
to the multiplication corresponding to
a path $f_1(h)=f_0+hp_1+\cdots$. Now we may apply the lemma to
$m^{(2)}_h$, and so on. On the $n$-st step we obtain
$m^{(n)}_h=(1+h^nD_n)...(1+h^2D_2)\circ m_h$ that
corresponds modulo $h^{n+1}$ to a path
$f_{n-1}(h)=f_0+hp_1+\cdots+h^{n-1}p_{n-1}+\cdots$.
Let $D=\lim (1+h^nD_n)...(1+h^2D_2)$, $f(h)=\lim f_{n-1}(h)$
in $h$-adic topology. It is clear that such limits exist.
We obtain that
$D\circ m_h=\mu_{f(h),h}$, which proves a).

The proof of b) using $(\Lambda^1\mm)^\lf=0$ is left to the reader.
\qed

\pat{\em Remark}\label{rem4.2}
Let $X_0$ be the variety of all nondegenerate Poisson brackets on $M$
(see \bref{ssec2.21}). As in the proof of \Prop{prop4.2} one can construct
a holomorphic family of $\Ug$ invariant associative multiplications of the form
\be{}\label{assom}
\mu_{p,t}(a,b)=ab+(t/2)p(a,b)+\sum_{n\geq 2}t^n\mu_{p,n}(a,b),\qquad p\in X_0.
\ee{}
The same argument as in the proof of \Prop{prop4.2a}
shows that such a family has the universal property: any 
$\Ug$ invariant deformation quantization on $M$ is equivalent
to the pullback of \bref{assom} by a unique formal path
in $X_0$.

\section{The $G$ invariant $\Phi_h$ associative quantization
in two parameters}

\pa\begin{propn}{prop4.3}Let $\g$ be a simple Lie algebra,
$M$ a semisimple orbit in $\g^*$. Let $v$ be the KKS Poisson bracket
on $M$. Let $f$ be an invariant bracket on $M$ satisfying
$\[f,f\]=-\ff_M$, $\[f,v\]=0$. Then  there exists a 
two parameter multiplication  $\mu_{t,h}$ on $\A$ 
\be{}\label{tpq}
\mu_{t,h}(a,b)=ab+(h/2)f(a,b)+(t/2)v(a,b)+
\sum_{k+l\geq 2} h^kt^l\mu_{k,l}(a,b)
\ee{}
which  is $\Ug$ invariant and
$\Phi_h$ associative.
\end{propn}

\begin{proof}
The existence of a multiplication which is $\Phi_h$ associative
up to and including $h^2$ terms is nearly identical to the 
proof of \Prop{prop4.2}.

 So, suppose we have a multiplication 
defined to order $n$,
$$\mu^{(n)}_{t,h}(a,b)=ab+h/2f(a,b)+t/2v(a,b)+
\sum_{2\leq k+l\leq n} h^kt^l\mu_{k,l}(a,b),$$
which is $\Ug$ invariant
and $\Phi_h$ associative to order $h^n$.
Consider the obstruction cochain,
$$obs_{n+1}=\sum_{k=0,\ldots, n+1} h^k t^{n+1-k}b_k.$$
The same argument as in the proof of \Prop{prop4.2} shows that
$obs_{n+1}$ is a Hochschild cocycle. This means that
all coefficients $b_k$ are Hochschild cocycles.
Hence, $b_k=\partial a_k+\beta_k$ for all $k$,
where $\beta_k\in (\Lambda^3\mm)^\lf$. Therefore,
$$obs_{n+1}=\partial a+\beta,$$
where $a=\sum h^kt^{n+1-k}a_k$,
$\beta=\sum h^kt^{n+1-k}\beta_k$.
The element $\beta$ is a cocycle from $\CC^3(n+1)$
(see Subsection \ref{sec3.3}). By \Prop{prop3.3} b)
there exists $\alpha\in\CC^2(n)$ such that 
$\[hf+tv,\alpha\]=\beta$. This shows that,
as in the proof of \Prop{prop4.2}, 
we can modify $\mu^{(n)}_{t,h}$ adding a multiple of $\alpha$ 
to get a new multiplication to order $n$ with $(n+1)$-st obstruction cocycle $obs_{n+1}$
being a Hochschild coboundary, $\partial a$. So, we are able
to continue the recursive construction of the
desired two parameter deformation.
\end{proof} 

\pat{\em Remark.} Let $\pi:Y\to X_0$ be the projection of the variety
of good brackets over $M$ to the variety of nondegenerate Poisson brackets
(see \Rem{rem2.1}).
In the similar way as \Prop{prop4.3}, one can prove the
existence of a family of multiplications of the form
\be{}
\mu_{f,t,h}(a,b)=ab+(h/2)f(a,b)+(t/2)(\pi f)(a,b)+
\sum_{k+l\geq 2}h^kt^l\mu_{f,k,l}(a,b), \quad f\in Y.
\ee{}
This family satisfies the universal property 
for two parameter quantizations, i.e., any two parameter
$\Ug$ invariant $\Phi_h$ associative quantization on $M$ 
of the form \bref{tpq} is the pullback of a two parameter formal
path in $Y$.                

\section{Polarization}

\pa We retain notations from the previous sections. Recall that $M=G/L$ 
where $G$ is a complex connected simple Lie group and
$L$ is a Levi subgroup. It is known that the natural projection
$\pi:G\to M$ is a holomorphic principal fiber bundle with
structure group $L$.  

The tangent space to $M$, $T(M)$, is the associated
vector bundle corresponding to the $\Ad$-action of $L$ on $\mm$, a
unique $\lf$ invariant subspace in $\g$ complement to $\lf$
(see Subsection \ref{ssec2.5}). According to \Prop{prop2.1},
$T(M)$ may be presented as a direct sum of subbundles,
$T(M)=\oplus_{\baa\in\Ob}T_\baa(M)$ where $T_\baa(M)$ is the
associated vector bundle corresponding to $\mm_{\baa}$.

Assigning to each $g\in G$ the horizontal subspace $g\mm$ provides
$G$ with an invariant connection $\nabla$. This connection defines
a $G$ invariant connection on any associated vector bundle over $M$.

\pa Let us choose $\Ob^+$, a system of positive quasiroots in $\Ob$.
Let the subbundle $T^+(M)=T_{\Ob^+}(M)$ ($T^-(M)=T_{\Ob^-}(M)$)
correspond to $\mm^+=\oplus_{\baa\in\Ob^+}\mm_\baa$ 
($\mm^-=\oplus_{\baa\in\Ob^-}\mm_\baa$).

By a realization $M$ as an orbit in $\g^*$ the decomposition
$T(M)=T_{\Ob^+}(M)\oplus T_{\Ob^-}(M)$ defines complement polarizations
of $M$ with respect to the KKS symplectic form on $M$.
These polarizations define two complement
foliations on $M$, which are fibrating with the natural projections
$M\to G/P^+$ and $M\to G/P^-$, respectively, where 
$P^+$, $P^-$ are upper and low parabolic subgroups containing $L$.

\pa\label{ssec6.2} Let $\A^+=\A_{\Ob^+}$ denote the sheaf of 
holomorphic functions $a$ on $M$ constant along the polarization
defined by $\Ob^+$, i.e.,
such that $\nabla_X a=0$ for any vector field $X\in T_{\Ob^+}(M)$.

\pa\begin{propn}{prop6.1} Let $\nu$ be a $\Ug$ invariant bidifferential
operator on $M$ vanishing on constants.
Let $\Ob^+$ be a system of positive quasiroots.
Then $\nu(a,b)=0$ for any sections $a,b\in\A_{\Ob^+}$.
\end{propn}

\begin{proof}
The connection $\nabla$ induces an equivariant isomorphism 
of $\A$ modules between the sheaf of differential operators on
$\A$ and the sheaf $ST^-(M)\ot_\A ST^+(M)$, where $ST^-(M)$ and $ST^+(M)$
denote the sheaves of symmetric tensors over $T^-(M)$ and $T^+(M)$.  
It provides an equivariant isomorphism between the space of invariant bidifferential
operators on $M$ and the space $((S\mm^-\ot_\C S\mm^+)^{\ot 2})^\lf$.
Thus one may regard $\nu$ as a sum of terms of the form
$A_1B_1\ot A_2B_2$ where $A_1,A_2\in S\mm^-$, $B_1,B_2\in S\mm^+$.
Since $\nu$ is invariant, $A_1B_1\ot A_2B_2$ must be of weight zero with
respect to the Cartan subalgebra $\hh\subset \lf$. Since $A_1B_1\ot A_2B_2$
is vanishing on constants, either $B_1$ or $B_2$ must belong to
a positive symmetric power of $\mm^+$, let such be $B_1$. But the
corresponding to $B_1$ differential operator takes the functions
of $\A^+$ to zero, therefore, the bidifferential operator
corresponding to  $A_1B_1\ot A_2B_2$ when applies  to the pair
$a,b\in \A^+$ gives zero, too. 
\end{proof}

\pa\begin{cor}{cor6.1}All  $\Ug$ invariant $\Phi_h$ associative
multiplications are trivial on $\A_{\Ob^+}$
for any choice of $\Ob^+$.
It means that for any of such a multiplication, $\mu$, one has
$\mu(a,b)=ab$ whenever $a,b\in \A_{\Ob^+}$.
\end{cor}

\begin{proof} Indeed, $\mu$ has the form 
$\mu(a,b)=ab+$\{bidifferential operators vanishing on constants\}.
So, the corollary follows from \Prop{prop6.1}.
\end{proof}

\pa Since $\mu$ from \Cor{cor6.1} when restricted to $\A^+$
coincides with the usual multiplication, it is associative in the usual
sense. On the other hand, $\mu$ is $\Phi_h$ associative, so the
usual and $\Phi_h$ associativities coincide on $\A^+$.
We will prove this fact independently in a more general setting.

Let $V$ be a representation of $L$. Denote by $V(M)$ the
corresponding associated vector bundle on $M$. 
When this does not lead to confusion 
we will use the same notation $V(M)$ for the sheaf of holomorphic
sections of the bundle $V(M)$.
Denote by $V^+(M)=V_{\Ob^+}(M)$   
the sheaf of holomorphic sections $v$ of $V(M)$ constant along the
polarization defined by $\Ob^+$, i.e.,
such that $\nabla_X v=0$ for any vector field $X\in T_{\Ob^+}(M)$.

It is clear that $V(M)$ is an $\A$ module under the natural multiplication 
$m: \A\ot V(M)\to V(M)$, $m(a,v)=av$, $a\in\A, v\in V(M)$. 
Since $\A$ is commutative, $V(M)$ may be considered as a two-sided module.
This multiplication is obviously associative, i.e.,
for sections $a,b\in\A$ and $v\in V(M)$ 
one has $(ab)v=a(bv)$ and $(av)b=a(vb)$.
The following proposition shows that this multiplication being restricted to
to $\A^+$ and $V^+(M)$ is also $\Phi_h$ associative.

\pa\begin{propn}{prop6.3}Let $V$  be a representation of $L$.
Let us choose a system of positive quasiroots $\Ob^+$.
Then for any sections $a,b\in\A^+$, $v\in V^+(M)$ one has
$abv=m\Phi_h(a\ot b\ot v)$ and $avb=m\Phi_h(a\ot v\ot b)$.
\end{propn}

\begin{proof} 
Let us prove the first relation. The second relations can be proven
similarly.
Let $\Omega^+$ be a system of positive roots for $\g$ which
projects on $\Ob^+$. Let $\frak{p}^+$ be the corresponding parabolic
subalgebra of $\lf$ and $\frak{u}^+$ the radical of $\frak{p}^+$, so
$\frak{p}^+=\lf\oplus\frak{u}^+$. Let $\frak{p}^-$ and $\frak{u}^-$
denote the corresponding opposite subalgebras, in particular,
$\g=\frak{p}^+\oplus\frak{u}^-$ and
$\g=\frak{p}^-\oplus\frak{u}^+$. 

Let $\pi:G\to M$ be the natural projection.
Let $U$ be an open set in $M$. Sections of $V(M)$ over $U$
can be identified with functions $f:\pi^{-1}(U)\to V$ such that
$f(gl)=l^{-1}f(g)$ for $l\in L$, $g\in G$. Sections of $V^+(M)$ must, 
in addition, satisfy the condition 
$E_\alpha f=0$, $\alpha\in\frak{u}^+$,
and the root vector $E_\alpha$ acts on $f$ as a complex
left-invariant vector field on $G$.
In particular, functions of $\A^+$ are identified with
functions $\psi$ over $\pi^{-1}$ such that $E_\alpha\psi=0$
for all $\alpha\in\frak{p}^+$.

Let us write $\Phi_h$ in the form
\be{}\label{Phiw}
\Phi_h=1\ot 1\ot 1+\sum_{k\geq 2} h^k\Phi_k^1\ot\Phi_k^2\ot\Phi_k^3, 
\ee{}
where
each $\Phi_k^i$, $i=1,2$, belongs to $S\frak{p}^+\ot S\frak{u}^-$
and $\Phi_k^3$ belongs to $S\frak{p}^-\ot S\frak{u}^+$
($S\frak{p}^+$ denote the space of symmetric tensors over $\frak{p}^+$,
and so on). The total degree of
each $\Phi_k^1\ot\Phi_k^2\ot\Phi_k^3$ is greater than zero.

Let us take $a,b\in\A^+$, $v\in V^+(M)$ and apply 
$\Phi_k^1\ot\Phi_k^2\ot\Phi_k^3$ to $a\ot b\ot v$
Suppose $\Phi_k^1(a)$ and $\Phi_k^2(b)$ are not equal to zero.
Then there are $x_1,x_2\in S\frak{u}^-$ such that 
$\Phi_k^1=A\ot x_1$, $\Phi_k^2=B\ot x_2$, where
$A,B\in S\frak{p}^+$. But since $\Phi_k^1\ot\Phi_k^2\ot\Phi_k^3$ is
an invariant element, it must be of degree zero under the Cartan subalgebra.
It follows that $\Phi_k^3$ has to be of the form 
$\Phi_k^3=C\ot x_3$ where $C\in S\frak{p}^-$, $x_3\in\frak{u}^+$
and $x_3\neq 0$. Hence, $\Phi_k^3(v)=0$.

So, we have proven that in the expression \bref{Phiw} all terms
except for the first when applying to $a\ot b\ot v$ are equal to zero.
\end{proof}

\section{The real case}

\pa Let $G$ be a real connected simple Lie group with
complexification $G^\C$. Let $\g_\R$ be the Lie algebra of $G$
with complexification $\g$.
Let $L$ be the stabilizer of a semisimple element $\lambda\in\g_\R^*$,
so that $M=G/L$ may be identified with the coadjoint orbit
passing through $\lambda$.
It is well known that $L$ is connected, therefore the complexification 
$L^\C\subset G^\C$ is meaningful.

Denote by $\lf_\R$ ($\lf$) the Lie algebra of $L$ ($L^\C$).
Note that $\lf$ is a Levi subalgebra in $\g$.
The natural embedding $M\to M^\C=G^\C/L^\C$ may be regarded as 
complexification of $M$.

Let $C^a(M)$ and $C^\infty(M)$ denote the sheaves of real analytic
and smooth complex valued functions on $M$.
The action of $G$ on $M$ defines a map of $\g_\R$ into
the Lie algebra of real vector fields on $M$, $\g_\R\to \Vect_\R(M)$,
that extends to a map $\g\to \Vect(M)$ of $\g$ into the Lie algebra of
complex vector fields on $M$.
It follows that $U(\g)$ acts on the sections of $C^a(M)$ and $C^\infty(M)$
as differential operators.  

As a consequence we get that all the $\Ug$ invariant $\Phi_h$
associative multiplications constructed in the previous sections
can be defined on the real manifold $M$.

\pa Let us choose a system of positive roots, $\Omega^+$, in $\g$.
Let $P$ be the corresponding parabolic subgroup of $G^\C$ with
Levi factor $L^\C$ and $\frak{p}$ its Lie algebra.
One has $\frak{p}=\lf\oplus\frak{u}^+$, where $\frak{u}^+$ is the 
nilradical of $\frak{p}$ 
We assume that  $\frak{p}$ is $\theta$-stable parabolic, i.e., 
satisfies the condition
\be{}\label{spar}
\g_\R\cap \frak{p}=\lf_\R.
\ee{}
Then the natural map 
\be{}\label{cstr}
M=G/L\to G^\C/P
\ee{}
is an inclusion and
the image is an open set. Thus the choice of $\Omega^+$ makes
$M$ into a complex manifold with holomorphic action of $G$.
The corresponding system of positive quasiroots, $\Ob^+$,
defines a complex polarization on $M$, whereas $\Ob^-$ defines
the complement polarization.

Note that for $\lf_\R$ a $\theta$-stable parabolic $\frak{p}$ exists if $\lf_\R$
is the centralizer of a semisimple element $x\in\g_\R$
such that $\ad(x)$ has imaginary eigenvalues. 
 
One can prove, \cite{Kn}, that the smooth functions on $M$ which are
constant along the polarizations defined by $\Ob^+$ and $\Ob^-$
coincide with holomorphic and antiholomorphic functions on $M$.

\pa Let $\tilde{S}$ denote the operator $F_h\sigma F_h^{-1}$,
where $F_h$ is from \Prop{propo2.1} and $\sigma$ is
the usual permutation, acting on the tensor product of
any two representations of $\Ug$. Let $B_h$ be a quantized algebra of  
functions on $M$. We say that the multiplication on $B_h$, $m_h$,
is $\tilde{S}$-commutative, if for any $a,b\in B$ one has
$m_h(a\ot b)=m_h\tilde{S}(a\ot b)$.  

\pa\begin{thm}{thm1}
Let $G$ be a real connected simple Lie group, $L$ a Lie subgroup 
which is a stabilizer of a semisimple element $\lambda\in\g_\R^*$,
and $M=G/L$.
Let $r\in\wedge^2\g$ be an $r$-matrix
and $\Uh$ the corresponding quantum group.
Let $X$ be the variety of $\ff$-brackets on $M$, as
in \Prop{prop4.2}, for $\ff=\[r,r\]$. 
Then there exists a universal family of multiplications on $C^\infty(M)$
of the form
\be{}
m_{f,h}(a,b)=ab+(h/2)(f-r_M)(a,b)+\sum_{n\geq 2}m_{f,n}(a,b), \quad f\in X,
\ee{}
which is $\Uh$ invariant and associative.  

Suppose the map \bref{cstr} induces a complex structure on $M$.
Then $m_{f,h}$ being restricted to the sheaf of holomorphic
functions is $\tilde{S}$ commutative. 
\end{thm}
  
\begin{proof} According to \Prop{propo2.2} we put
$m_{f,h}=\mu_{f,h}F^{-1}_h$, where the multiplication
$\mu_{f,h}$ is from \Prop{prop4.2}. The functions of $C^a(M)$
are restrictions to $M$ of holomorphic functions
on $M^\C$. It follows that $m_{f,h}$ is a well defined
multiplication on $C^a(M)$. Since the coefficients by $h^n$ 
in $m_{f,h}$ are bidifferential operators, this multiplication
is defined, actually, for smooth complex valued functions on $M$.
The universality of $m_{f,h}$ follows from the universality
of $\mu_{f,h}$, see \Prop{prop4.2a}.
The $\tilde{S}$-commutativity of $m_{f,h}$ for holomorphic
functions follows directly from the commutativity of $\mu_{f,h}$
for such functions, see \Cor{cor6.1}.
\end{proof} 

\pa As a consequence, we obtain the statement reverse to \Cor{coro2.3}:
any Poisson bracket on $M$ of the form \bref{dopq} admits a
$\Uh$ invariant quantization. The proof is analogous to \Cor{cor4.2}.

The following theorem contains, in particular, the statement reverse
to \Cor{coro2.4}.

\pa\begin{thm}{thm2}
Let $M$ be as in \Thm{thm1}.
Let $v$ be the KKS Poisson bracket on $M$.
Let $r\in\wedge^2\g$ be an $r$-matrix
and $\Uh$ the corresponding quantum group. 
Let $p$ be a Poisson
bracket on $M$ of the form $p=f-r_M$, where $f$ satisfies
$\[f,f\]=\[r,r\]_M$ and $\[f,r_M\]=0$.
Then there exists a two parameter $\Uh$ invariant associative 
multiplication on $C^\infty(M)$ of the form
\be{*}
m_{t,h}(a,b)=ab+h/2(f-r_M)(a,b)+t/2v(a,b)+\sum_{k+l\geq 2}h^kt^lm_{k,l}(a,b).
\ee{*}
Suppose the map \bref{cstr} induces a complex structure on $M$.
Then $m_{f,h}$ being restricted to the sheaf of holomorphic
functions is $\tilde{S}$ commutative. 
\end{thm}

\begin{proof} We put $m_{t,h}=\mu_{t,h}F_h^{-1}$ and use the argument
as in the proof of \Thm{thm1}.
\end{proof}

\pat{\em Remark.}\label{rem8.4} As follows from \Cor{cor6.1}, any $\Uh$ invariant
multiplications, in particular the multiplications
from Theorems \ref{thm1} and \ref{thm2}, being restricted to
holomorphic functions are equal to $m_0F_h^{-1}$, where $m_0$ is
the usual multiplication.

\section{The quantization of vector bundles}
\pa 
Let $\rho:L\to GL(V)$ be a representation of $L$ in a complex
vector space $V$. Then $\rho$ extends holomorphically to $L^\C$.
Let $\rho_\C$ denote such an extension.
The vector bundle on $M$, $V(M)$, associated with $\rho$
is the restriction of the vector bundle on $M^\C$, $V(M^\C)$, associated with
$\rho_\C$. Let us choose a system of positive quasiroots, $\Ob^+$,
and the corresponding parabolic subgroup, $P\supset L^\C$.
Let us assume that the map \bref{cstr} defines the complex
structure on $M$. 
Let $\rho_P$ denote the extension of $\rho_\C$
to $P$, which is trivial on the unipotent radical of $P$.
Then $V(M)$ is the pullback of the vector bundle on $G^\C/P$
associated by $\rho_P$. 
So, $V(M)$ acquires the structure of a holomorphic vector bundles on $M$.
The holomorphic sections of $V(M)$ form a sheaf $V^+(M)$ whose
sections are the restrictions to $M$
of sections of $V^+(M^\C)$ constant along the polarization
defined by $\Ob^+$. 

In the following we fix a system $\Ob^+$ and the corresponding
complex structure on $M$. Given a representation $V$ of $G$,
we denote by $V(M)$ and $V^+(M)$ the sheaves of smooth and
holomorphic sections of the associated to $V$ vector
bundle on $M$.

\pa\begin{defn}{def8.1} Let $\A\subset C^\infty(M)$ be a subsheaf
of algebras
and $E$ a sheaf of $\A$ modules with a $\Ug$ action.
Let $m_0:\A\ot\A\to\A$ and $n_0:\A\ot E\to E$ 
denote the multiplication in $\A$ and the action of $\A$ on $E$.
Let $\A_h$ be a $\Uh$ invariant quantization of $\A$ with
the deformed multiplication $m_h=m_0+h m_1+o(h)$.
We say that $E[[h]]$ is a quantizqtion of $E$ as a 
$\A_h$ module if a deformed $\Uh$ invariant action
$n_h=n_0+h n_1+o(h):\A\ot E\to E[[h]]$ is given, which makes $E[[h]]$ into
an $\A_h$ module. In particular, it means that the associativity holds:
$n_h(m_h(a,b),x)=n_h(a,n_h(b,x))$ for $a,b\in\A$, $x\in E$.

Since $\A$ is commutative, $E$ is, in fact, a two-sided module.
So, in the similar way one defines a quantization of $E$ as
a two-sided $\A_h$ module.
\end{defn}

\pat{Example.}\label{exa8.5} Let $\A^+$ be the sheaf of holomorphic functions on $M$.
Then, as follows from \Rem{rem8.4}, there exists a unique $\Uh$ invariant
quantization of $\A^+$, $\A^+_h$, and it has the multiplication of the form
$m_h=m_0F_h^{-1}$. Let $V$ be a representation of $L$.
Then \Prop{prop6.3} and the argument of \Thm{thm1} show that 
the sheaf of holomorphic sections $V^+(M)$ can be uniquely quantized 
as a left and even as a two-sided module. We denote this quantization by $V^+_h(M)$.
The left (right) multiplication by elements of $\A$ has the form
$n_h(a\ot x)=n_0F_h^{-1}(a\ot x)$ ($n_h(x\ot a)=n_0F_h^{-1}(x\ot a)$)
for $a\in\A$, $x\in E$.

The following proposition shows that the sheaf of smooth
sections $V(M)$ admits a $\Uh$ invariant quantization.

\pa\begin{thm}{thm3}Let $\A=C^\infty(M)$ and $\A_h$ be
a $\Uh$ invariant quantization of $\A$. Let $V$ be a representation of $L$.
Then there exists a $\Uh$ invariant quantization of $V(M)$ as
a left $\A_h$ module.
\end{thm} 

\begin{proof} We have $V(M)=\A\ot_{\A^+}V^+(M)$. Let
$V_h(M)=\A_h\ot_{\A^+_h}V^+_h(M)$, where $\A_h$ is considered as 
a right and $V^+_h(M)$ as a left $\A^+_h$ module
(see \Exa{exa8.5}). It is clear that $V_h(M)$ is the required quantization.
\end{proof}   

\pat{The two-sided quantization.} In general, it is not clear whether a 
two-sided quantization of $V(M)$ exists. However we will show that
there is a quantization, $\A_h$,  of the sheaf of smooth functions on $M$
such that for any representation $V$ there exist a quantization of $V(M)$ 
as a two-sided $\A_h$ module.

Let us construct the quantization $\A_h$.  
Let 
$$R=F_h{\rm e}^{h{\mathbf t}/2}F_h^{-1}=
R^\prime_i\ot R^{\prime\prime}_i\in \Uh\ot \Uh,$$ 
where ${\mathbf t}\in\g\ot\g$ is the split Casimir, be the R-matrix 
(summation by $i$ is assumed). 
It satisfies the property, \cite{Dr2},
\be{*} 
\Delta^\prime(x)=R\Delta(x)R^{-1}, \quad x\in \Uh,
\ee{*}
where $\Delta$ is the comultiplication in $\Uh$ and $\Delta^\prime$  
the opposite one,
\be{*} 
(\Delta\ot 1)R&=&R^{13}R^{23}=R^\prime_i\ot R^\prime_j\ot R^{\prime\prime}_i
R^{\prime\prime}_j \\
(1\ot\Delta)R&=&R^{13}R^{12}=R^{\prime}_iR^{\prime}_j\ot R^{\prime\prime}_j
\ot R^{\prime\prime}_i, 
\ee{*}
and
\be{*} \label{coun}
(1\ot \varepsilon)R=(\varepsilon\ot 1)R=1\ot 1, 
\ee{*}
where $\varepsilon$ is the counit in $\Uh$.

The element $R$ defines the $\Uh$ equivariant map $S:E\ot F\to F\ot E$, 
$a\ot b\mapsto \sigma R(a\ot b)$, $\sigma$ is the usual permutation,
for any $\Uh$ modules $E$ and $F$.

Let $\A'_h=\A_h^+\ot_\C \A^-_h$, where $\A_h^-$ is a unique
$\Uh$ invariant quantization of the sheaf $\A^-$ of antiholomorphic functions on $M$.
We provide $\A'_h$ with the structure of a sheaf of algebras in the
following way. For $a=a_1\ot b_1, b=a_2\ot b_2\in \A'_h$, we put
$m_h(a,b)=a_1a'_2\ot b'_1b_2$, where $a'_2\ot b'_1=S(b_1\ot a_2)$
and $a_1a'_2$ and $b'_1b_2$  means the multiplications in $\A_h^+$
and $\A_h^-$, respectively. It easily follows from the above properties of $R$ that the
multiplication $m_h$ is $\Uh$ invariant and associative on $\A'_h$ and is
presented as a power series in $h$ with coefficients being
bidifferential operators on $\A$. Since bidifferential operators on smooth
functions are fully 
defined by their values on $\A^+\ot\A^-$, 
this multiplication can be extended to the whole algebra $\A$
of smooth functions on $M$.

One can show that the Poisson bracket of the obtained quantization $\A_h$
is the bracket reduced to $M$ from the Poisson bracket $r'-r''$ on the group $G$,
see \Rem{rem2.4}.

\pa\begin{thm}{thm4} Let $\A_h$ be the quantization of the sheaf of
smooth functions on $M$ constructed above. Let $V$ be a representation
of $L$. Then there exists a quantization of the sheaf of smooth sections
of $V(M)$ as a two-sided $\A_h$ module.
\end{thm}

\begin{proof} Let $V'_h(M)=V^+_h(M)\ot_\C\A^-_h$. Let us define left
and right multiplications of elements of $V'_h(M)$ by elements
of $\A'_h=\A_h^+\ot_\C \A^-_h$.
Let $a=a_1\ot b_1\in\A'_h$ and $x=x_1\ot b_2\in V'_h(M)$.
Put $n^{left}_h(a\ot x)=a_1x'_1\ot b'_1b_2$, where
$x'_1\ot b'_1=S(b_1\ot x_1)$,
and $n^{right}_h(x\ot a)=x_1a'_1\ot b'_2b_1$, where
$a'_1\ot b'_2=S(b_2\ot a_1)$.
Here $a_1x'_1$ means, for example, the multiplication
when $V^+_h(M)$ is considered as a left $\A^+_h$ module.   
It easily follows from the above properties of $R$ that the multiplications
$n^{left}_h$ and $n^{right}_h$ make  $V'_h(M)$ into
a two-sided $\A'_h$ module. The same argument as in the
proof of \Thm{thm3} shows that those multiplications
define, in fact, the structure of a $\A_h$ module on $V(M)[[h]]$.  
\end{proof}

\end{document}